\documentclass[final]{siamart1116}
\usepackage{indentfirst,latexsym,bm,amssymb,amsfonts,ifthen,color}
\usepackage{subfigure}
\usepackage{graphicx}
\usepackage{algorithm,algorithmic}
\usepackage{amsmath}
\usepackage{mathrsfs}
\usepackage{enumerate,enumitem}
\usepackage{multirow,diagbox}

\newtheorem{remark}{\bf Remark}

\newcommand{\N}{\mathcal{N}}
\newcommand{\T}{\mathcal{T}}

\def\O{\mathscr{O}}
\def\R{\mathcal{R}}
\def\M{\mathcal{M}}
\def\RR{\mathbb{R}}
\def\E{\mathcal{E}}
\def\RO{\mathcal{R}_\Omega}
\def\DO{D_\Omega}
\def\TDO{\tilde{D}_\Omega}
\def\Tr{\text{Trace}}

\title{Solving Partial Differential Equations on Manifolds From Incomplete Inter-Point Distance
}

\author{Rongjie Lai
\thanks{Department of Mathematics, Rensselaer Polytechnic Institute, Troy, NY 12180,
         U.S.A. ({\tt lair@rpi.edu}). The research of Rongjie Lai is partially supported by NSF grant DMS--1522645.}
        \and Jia Li
        \thanks{Department of Mathematics, Rensselaer Polytechnic Institute, Troy, NY 12180,  ({\tt lij25@rpi.edu}).}
        }
\date{}

\begin{document}
\maketitle

\begin{abstract}
Solutions of partial differential equations (PDEs) on manifolds have provided important applications in different fields in science and engineering. Existing methods are majorly based on discretization of manifolds as implicit functions, triangle meshes, or point clouds, where the manifold structure is approximated by either zero level set of an implicit function or a set of points. In many applications, manifolds might be only provided as an inter-point distance matrix with possible missing values. This paper discusses a framework to discretize PDEs on manifolds represented as incomplete inter-point distance information. Without conducting a time-consuming global coordinates reconstruction, we propose a more efficient strategy by discretizing differential operators only based on point-wisely local reconstruction. Our local reconstruction model is based on the recent advances of low-rank matrix completion theory, where only a very small random portion of distance information is required. This method enables us to conduct analyses of incomplete distance data using solutions of special designed PDEs such as the Laplace-Beltrami (LB) eigen-system. As an application, we demonstrate a new way of manifold reconstruction from an incomplete distance by stitching patches using the spectrum of the LB operator. Intensive numerical experiments demonstrate the effectiveness of the proposed methods.
\end{abstract}

\begin{keywords}
Manifolds, Laplace-Beltrami eigenproblem, Eikonal equation, Low-rank matrix completion.
\end{keywords}

\begin{AMS}
  65D18, 65D25, 65N25
\end{AMS}

\section{Introduction}
\label{sec:intro}
With the rapid development of advanced data acquisition technology, processing and analyzing data sampled on 3D shapes or even higher dimensional geometric objects becomes ubiquitous tasks such as those used in a 3D camera, medical imaging, protein structuring, social network analysis and many others~\cite{faugeras1986representation,axelsson1999processing,waterman1995introduction,starck1998image,crippen1988distance,berger1999reconstructing,scott2011sage,ji2004sensor,biswas2006semidefinite}.
PDE and variational PDE based methods have made great success to handle problems in signal and image processing which can be viewed as data on Euclidean domains. It is natural to consider PDE based methods to analyze and process signals on a general manifold and to understand geometric structures hidden in the data. Besides the classical implicit methods~\cite{Osher:88,Bertalmio:2000, Bertalmio:2002}, finite difference methods~\cite{pinkall1993computing,taubin2000geometric,Meyer:2003,xu2004convergent}, finite element methods~\cite{Reuter:06,lai2011framework,dziuk2013finite} and parameterization methods~\cite{stam2003flows,spira2007geometric,wang2007brain,lui2008variational} for solving differential equations on surfaces in $\RR^3$, there has been increasing interests of solving PDEs on general $d$-dimensional manifold in $\RR^p$ and their applications to data analysis. For instance, a diffusion geometry framework is developed to investigate the geometric structure of data based on solving Laplace-Beltrami (LB) eigenproblem using integral kernel methods~\cite{Belkin:ML2004,Belkin:09clp,coifman2006diffusion}. More recently, a moving least square method and a local mesh method are considered to intrinsically solve different types of PDEs on manifolds represented as point clouds and its applications to the geometric understanding of point cloud data~\cite{lai2013local,Liang:CVPR2012,liang2013solving,lai2017multi}.

All the aforementioned methods of solving PDEs on a general manifold $\M \subset \RR^p$ are commonly considered $\M$ is sampled as a set of points $\{\bm{x}_i\in\RR^p\}_{i=1}^n$, referred as a {\it point cloud}, and discretization of differential operators or approximation of integral equations are relied on available coordinates information of $\{\bm{x}_i\}_{i=1}^n$. However, there are many applications whose input information has no point coordinates but only an incomplete inter-point distance $(d(\bm{x}_i,\bm{x}_j))$. Examples include incomplete distance information from sensor network localization~\cite{ji2004sensor,biswas2006semidefinite}, protein structuring from NMR spectroscopy \cite{crippen1988distance,berger1999reconstructing} and global position from local distance configuration of cities~\cite{singer2008remark}. A well-known distance geometry problem~\cite{crippen1988distance,mucherino2012distance} is to find the global configuration of data based on the incomplete distance information. The objective in this paper is different from the canonical distance geometry problem. We would like to develop numerical methods for solving PDEs on manifolds represented as incomplete inter-point distance. One natural way to approach this problem is first to apply a global reconstruction algorithm to obtain a point cloud representation of the input manifold, then established numerical methods for PDEs on point clouds can be directly applied. However, the global coordinate reconstruction might be very time consuming as it involves with semi-positive definite programming whose size depends on the number of points and could be very large in practice~\cite{biswas2006semidefinite}.

In this paper, we propose a different strategy to solve PDEs on manifolds represented as incomplete inter-point distance without conducting global coordinates reconstruction. Our idea is based on two intrinsic properties of differential operators on manifolds. Namely, the definition of a differential operator is only point-wisely depending on local information of the manifold and is invariant under different choices of local coordinates. This motivates us to only consider to conduct point-wisely local coordinate reconstruction for the associated local neighborhood, then we can point-wisely approximate differential operators based on the reconstructed local coordinates. After that, numerical solver for differential equations can be constructed.

Inspired by the classical multi-dimensional scaling~\cite{kruskal1978multidimensional}, the full distance matrix is one-to-one corresponding to the Gram matrix which can be further used to determine coordinates by its eigen-decomposition. More importantly, the Gram matrix can be essentially viewed as an inner product matrix after certain coordinates shift. As long as the matrix size, namely the size of a given local neighborhood, is larger than the embedding dimension of the local neighborhood, then the corresponding Gram matrix has to be a low-rank matrix. This suggests us to use low rank as a prior knowledge to reconstruct the Gram matrix based on constraints of available distance information. Adapted from the recent advances of the low-rank matrix completion theory~\cite{CandesRe2008}, we consider a nuclear norm regularized convex optimization problem to reconstruct local coordinates based on available distance information. Once local coordinates can be obtained, we apply the intrinsic methods of approximating differential operators on point clouds developed in~\cite{lai2013local,Liang:CVPR2012,liang2013solving} to discretize the desired differential equation. These intrinsic methods can be used to discretize different types of PDEs including parabolic, elliptic and hyperbolic PDEs. Our method can be viewed as natural extensions of these two methods to a new data structure where no coordinate information but only partial inter-point distance information of point clouds is provided. The outline of our strategy can be  as summarized as follows:

\begin{description}
\item[Step 0] For the $i$-th point, chose its K-nearest neighborhood (KNN) $N(i)$ based on the given incomplete distance or from the prior information if its available.

\item[Step 1] Applying the matrix completion method discussed in Section~\ref{sec:CoordinateRecon} to reconstruct local coordinates for KNN of the $i$-th point.

\item[Step 2]  Applying the MLS method or the local mesh method in~\cite{lai2013local,Liang:CVPR2012,liang2013solving} to approximate the desired differential operators at the $i$-th point. This provides the $i$-th row of the discretized matrix representation of the desired differential equation.

\end{description}

An immediate advantage of this approach is to save computation time by avoiding the global coordinate reconstruction, which can reduce the complexity quadratically to linearly scaling with the total number of points. More details about this method will be discussed in Section~\ref{sec:SolvePDEs} and will be verified in our numerical experiments. Furthermore, this approach enables us to conduct geometric understanding of data without global coordinates reconstruction. Examples include global pattern extraction, comparison and classification as many existing methods conducted using results of differential equations~\cite{Shi:08a,Lai:2010CVPR,lai2017multi}. As a byproduct of global information from PDEs, we also propose a new method for reconstructing manifolds by stitching its local patches. This new method is much more efficient than the way of direct global reconstruction using matrix completion method. Moreover, it can also overcome possible reconstruction failure using global matrix completion due to the coherent missing information.

The rest of this paper is organized as follows. In section~\ref{sec:CoordinateRecon}, we propose a low-rank matrix completion model to reconstruct local coordinates and design an algorithm to solve the proposed convex optimization problem based on operator splitting and the alternating direction method. After that, section~\ref{sec:SolvePDEs} is devoted to discuss point-wisely approximating differential operators based on the intrinsic methods proposed in~\cite{lai2013local,Liang:CVPR2012,liang2013solving}. We also discuss our two model PDEs, the Laplace-Beltrami (LB) eigenvalue problem and the Eiknoal  equation for solve geodesic distance on manifolds. As an application of global information using solutions of LB eigenvalue problem, we propose a new method for reconstructing manifolds by stitching local patches in the LB frequency space in section~\ref{sec:MfdStitching}. All numerical experiments about the effectiveness of our local coordinates reconstruction model, the accuracy and robustness of the PDEs solvers, and demonstration of manifold reconstruction are discussed in section~\ref{sec:experiments}. Finally, we conclude our work in section~\ref{sec:conclusion}.

\section{Local coordinates reconstruction via matrix completion}
\label{sec:CoordinateRecon}
In this section, we first review the concept of the classical multidimensional scaling (MDS). Inspired by the classical MDS, we propose a method of reconstructing local coordinates from the given incomplete inter-point distance matrix using the low-rank matrix completion~\cite{CandesRe2008}. A numerical algorithm based on the augmented Lagrangian is also designed to solve the proposed convex optimization problem. As we described before, our strategy is to conduct local coordinates reconstruction for the KNN $N(i)$ of the $i$-th point. For convenience, we write $\{\bm{x}_1,\bm{x}_2,\ldots \bm{x}_\ell\}$ as $\ell$ points for $N(i)$ as the following model does not rely on a typical choice of $N(i)$.

\subsection{Classical multidimensional scaling}
\label{subsec:MDS}
The classical MDS is to find a configuration of a set of points $\{\bm{x}_1,\bm{x}_2,\ldots \bm{x}_\ell\}$ in $\RR^p$ from the given distance matrix $D = (d_{ij}^2)_{\ell\times \ell}$ \footnote{For convenience, each entry of $D$ is specified by the square of distance.} such that the Euclidean distance matrix $\left(\|\bm{x_i} - \bm{x}_j\|^2_2 \right)_{ij}$ among these points is as close as possible to the given distance matrix $D$~\cite{kruskal1978multidimensional,borg2005modern}.

\begin{definition}
\label{def:Eucli}
A distance matrix $D = (D_{ij})_{\ell\times \ell}$ is called Euclidean if there exist some points $\bm{x}_1,\bm{x}_2,\ldots,\bm{x}_\ell \in \mathbb{R}^p$ such that $D_{ij}=(\bm{x}_i-\bm{x}_j)^{\top}(\bm{x}_i-\bm{x}_j)$. Let's write $X=(\bm{x}_1,\bm{x}_2,\ldots \bm{x}_\ell)^{\top}$ as an $\ell \times p$ matrix and denote the centering matrix by $\displaystyle H=I_l-\frac{1}{l} \bm{1} \bm{1}^{\top}$.
A Gram matrix $B = (B_{ij}) $ is defined as $B= -\frac{1}{2}HDH$, equivalently, $D_{ij} = B_{ii} + B_{jj} - B_{ij} - B_{ji}$.
\end{definition}

Note that the Euclidean distance matrix is translation invariant, thus, the matrix $H$ essentially translates the point set centering at the origin. It is well-known that $D$ is Euclidean if and only if the Gram matrix $B$ is positive semidefinite. Furthermore, since $B = HXX^{\top} H$, it is clear that $\text{rank}(B)=\text{rank}(XX^{\top})=\text{rank}(X)=p$. In our case, we consider the number of points is much larger the dimension of points. This suggests the low-rank property of $B$ as its prior knowledge.
If $D$ is Euclidean, and thus $B$ is positive semidefinite, then the eigen-decomposition of $B = U \Lambda U^{\top}$ provides the coordinate reconstruction $X = U_p\Lambda^{1/2}_p$, where $\Lambda_p$ is a diagonal matrix formed by the largest $p$ eigenvalues of $X$ and $U_p$ are the corresponding eigenfunctions.

\subsection{Local coordinates reconstruction via positive-semidefinite matrix completion}
Different from the MDS described as the above, one of the major challenges in our assumption is that only a small portion of $D$ is available. Generally speaking, it is impossible to recover the whole matrix only from its small portion of information. We here assume the information $D$ is randomly missing. Different models of semi-definite programming have been considered~\cite{alfakih1999solving,biswas2006semidefinite} in this scenario. Here, we consider a different model inspired by the low-rank prior information of the Gram matrix and the recent advances from the low-rank matrix completion theory~\cite{CandesRe2008}. In other words, we are seeking for a symmetric positive-semidefinite low-rank matrix $B$ satisfying $B_{ii}+B_{jj}-B_{ij} - B_{ji} = D_{ij}$ for $(i,j)$ in the available portion. Mathematically,  let's write $\Omega\subset \{(i,j) ~|~  1\leq j< i < \ell \}$ as the available index set and write $\mathbb{S} = \{B = (b_{ij}) \in\RR^{\ell\times \ell}~|~ B = B^T\}$.
We further write the available part of $D$ as a vector $D_\Omega \in \RR^{|\Omega|}$ and define a restriction operator $\mathcal{R}_{\Omega}:\mathbb{S} \rightarrow \RR^{|\Omega|} $ by $\mathcal{R}_{\Omega}(B)_{ij} = B_{ii}+B_{jj}-B_{ij} - B_{ji}$ with $(i,j)\in\Omega$. Our problem becomes to find a positive semidefinite matrix $B\in\mathbb{S}$ satisfying $\mathcal{R}_{\Omega}(B) = D_\Omega$.
It is clear that this problem is underdetermined and there are infinitely many solutions, most of which are far from the ground truth of $B$ and thus not applicable. Regularization method is a common way to handle this type of situation. Based on the prior knowledge of $B$ is a low rank matrix without knowing the exact rank $p$, we here consider the following optimization problem.
\begin{equation}
\label{eqn:ReconModel_B1}
\min_{B\in\mathbb{S}}  \mathrm{Rank}(B), \quad\text{s.t.} \quad B \succeq 0  \quad \&\quad   \mathcal{R}_\Omega(B) = D_\Omega
\end{equation}
Inspired by the matrix completion theory~\cite{CandesRe2008}, we relax the above NP-hard problem by replacing $\mathrm{Rank}(B)$ as the nuclear norm $\|B\|_*$. In addition, the input distance matrix is invariant under translation of the reconstructed coordinates $\{\bm{x}_i\}$. This translation ambiguity can be fixed by requiring $\sum_i \bm{x}_i=\bm{0}$. This condition is equivalent to row sum of $B$ is zero, which can also be obtained from the definition of the Gram matrix $B$. Therefore, we propose the following matrix completion model:
\begin{equation}
\label{eqn:ReconModel_B2}
\min_{B\in\mathbb{S}}  \|B\|_*, \quad\text{s.t.} \quad B \succeq 0,   ~~  \mathcal{R}_\Omega(B) = D_\Omega ~~ \&~~ B \bm{1} = 0
\end{equation}
where $\bm{1}$ denotes a column vector with all elements constant $1$. Once the above Gram matrix is obtained, coordinates can be reconstructed by the eigen-decomposition of $B$ as the method used in the classical MDS described in section~\ref{subsec:MDS}. Simultaneously, the eigen-decomposition of the reconstructed $B$ clearly shows the exact dimension $p$ by counting the cardinality of non-zero singular value. We would like to remark that the constraints and the restriction operator are different from the setting in~\cite{CandesRe2008} and it is highly nontrivial to check the restricted isometry property~\cite{recht2010guaranteed}, therefore, theoretical analysis conducted in~\cite{CandesRe2008,recht2010guaranteed} can not be directly applied in our problem. By considering a dual basis method, we can also theoretically show the exact completion can be achieved under certain coherence condition. 
More details about the theoretical analysis of the model will appear in our ongoing work~\cite{lai2017Exact}. Additionally, in a related work~\cite{weinberger2004learning}, the authors propose a semidefinite embedding (SDE) model as a kernel learning method for nonlinear dimension reduction, where no random missing but only KNN information of $D$ is considered, and the rank minimization is not considered but maximization of $Trace(B)$ is proposed for maximizing the variance in feature space.

\subsection{Numerical Algorithm}
For a convenience of designing a numerical algorithm of model \eqref{eqn:ReconModel_B2}, we first introduce some notations. Note that the set of symmetric matrix $\mathbb{S}$ can be viewed as $\RR^{\ell(\ell+1)/2}$
due to a natural isomorphism $\zeta:\mathbb{S} \rightarrow \RR^{\ell(\ell+1)/2}, B\mapsto \hat{B} = (b_{11},\cdots,b_{\ell1},b_{22},\cdots,b_{\ell2},\cdots, b_{\ell\ell})^T$. We denote $\iota = \zeta^{-1}$ and further merge the two linear constraints by defining the following linear operator:
\begin{equation}
\label{eqn:A}
\mathcal{A} : \RR^{l(l+1)/2}  \rightarrow \mathbb{R}^\Omega \times \mathbb{R}^\ell, \qquad \hat{B} \mapsto  \left ( \RO\circ \iota(\hat{B}), \iota(\hat{B}) \bm{1} \right)
\end{equation}
and extend the vector $\DO$ as $\TDO = (D_\Omega, \vec{0})$. Note that $\|B\|_*$ is the same as $\Tr(B)$ as $B$ is positive semi-definite. Therefore, the proposed model \eqref{eqn:ReconModel_B2} can be written as
\begin{equation}
\label{eqn:ReconModel_B3}
\min\limits_{\hat{B}\in\RR^{k(k+1)/2}}\Tr(\iota(\hat{B})), \quad
 \text{s.t.} \quad \mathcal{A}\hat{B} = \TDO,  \quad
\iota(\hat{B}) \succeq 0
\end{equation}
The model \eqref{eqn:ReconModel_B3} is a semi-definite programming (SDP) problem. It is well-known that the SDP problem can be solved by interior-points method \cite{SDO,SDP} as a canonical choice. However, when the number of constraint is with the order $\O(\ell^2)$ as in \eqref{eqn:ReconModel_B3}, the computational time can be with the complexity $\O(\ell^6)$~\cite{AltSDP-WenGoldfarbYin-2009}. Thus, we use the alternating direction method \cite{AltSDP-WenGoldfarbYin-2009} to save both the time and memory consumption. We first introduce an auxiliary variable $\hat{C}=\hat{B}$ and write \eqref{eqn:ReconModel_B3} as: 
\begin{equation}
\label{eqn:ReconModel_AL}
\begin{aligned}
\min\limits_{\hat{B},\hat{C}\in\RR^{n(n+1)/2}}\Tr(\iota(\hat{B})),
 \quad \text{s.t.} \quad \mathcal{A}\hat{B} = \TDO, \quad
 \iota(\hat{C}) \succeq 0, \quad
 \hat{B}=\hat{C}
 \end{aligned}
\end{equation}
Then we introduce an augmented Lagrangian of model \eqref{eqn:ReconModel_AL} as follows:
\begin{equation}\label{eqn:model global aug lag}
\mathcal{L}_{\mu_1,\mu_2}(\hat{B},\hat{C}; H_1,H_2)= \Tr(\iota(\hat{B})) + \frac{\mu_1}{2}\| \mathcal{A}\hat{B}-\TDO+ H_1 \|_2^2+ \frac{\mu_2}{2}\| \hat{B} - \hat{C} + H_2 \|_2^2.
\end{equation}
where $H_1$ and $H_2$ are two dual variables. The variable splitting enables us to alternatively optimize $\hat{B}$ and $\hat{C}$. As \eqref{eqn:ReconModel_AL} is convex, the saddle point of $\mathcal{L}_{\mu_1,\mu_2}(\hat{B},\hat{C}; H_1,H_2)$ is the solution which can be obtained by the following iterative procedure:
\begin{equation}
\label{ReconModelAL_ alg}
\begin{cases}
 \displaystyle \hat{B}^{k+1}=\arg\min_{\hat{B}}\Tr(\iota(\hat{B})) + \frac{\mu_1}{2}\| \mathcal{A}\hat{B}-\TDO+ H_1 \|_2^2+ \frac{\mu_2}{2}\| \hat{B} - \hat{C} + H_2 \|_2^2 \\
 \displaystyle  \hat{C}^{k+1}=\arg\min_{\hat{C}} \frac{\mu_2}{2}\|\hat{B}^{k+1}- \hat{C} + H_2^{k} \|^2 ,\quad \iota(\hat{C}) \succeq 0,\\
H_1^{k+1}=H_1^{k}+ (\mathcal{A} \hat{B}-\TDO), \\
H_2^{k+1}=H_2^{k}+  \hat{B} - \hat{C}, \\
\end{cases}
\end{equation}
Note that the model \eqref{eqn:ReconModel_AL} is convex, therefore convergence of the above method can be guaranteed~\cite{tseng2001convergence, AltSDP-WenGoldfarbYin-2009}.

The first subproblem is convex and differentiable. From its first order optimality condition, $\hat{B}^{k+1}$ solves the following linear system:
\begin{equation}\label{linearsys1step}
 (\mu_1 \mathcal{A}^*\mathcal{A}+\mu_2) \hat{B}^{k+1}= \mu_1 \mathcal{A}^*(\TDO-H_1^{k}) +\mu_2(\hat{C}^k-H_2^{k})- \zeta(I_n) .
\end{equation}
where $\mathcal{A}^*$ is the adjoint operator of $\mathcal{A}$ defined as follows. It is clear that for any given $\tilde{U} =  (U_\Omega, \vec{q}) \in \mathbb{R}^\Omega \times \mathbb{R}^\ell$ with $U_\Omega= \{u_{ij}\}_{(i,j)\in\Omega} $, the conjugate of $\mathcal{A}$ is given by $\mathcal{A}^* \tilde{U}= \zeta(U) + \zeta(Q)$. Here we define symmetric matrices $U$ and $Q$ by:
\begin{equation}\label{eqn:A*}
 U_{ij} = \begin{cases}
\displaystyle \sum_{(i,k)\in\Omega} u_{ik}+ \sum_{(k,i)\in\Omega}u_{ki}, & \text{if} \  i = j\\
  - 2 u_{ij}, &\text{if} \ i > j \end{cases}  \quad \text{and} \quad
Q_{ij} = \begin{cases}
 q_i, & \text{if} \  i = j\\
 q_i + q_j, &\text{if} \ i > j \end{cases}
 \end{equation}
and $U_{ij} = 0 $ for any $(i,j) \notin \Omega$.
As $\mu_1 \mathcal{A}^*\mathcal{A}+\mu_2 I_\ell$ is a symmetric and positive definite matrix, there are numerous numerical solvers to approach the solution of the above problem. Here, we typically chose the conjugate gradient method to solve \eqref{linearsys1step}.

The second subproblem of \eqref{ReconModelAL_ alg} is to find the nearest point of $\hat{B}^{k+1}+H_2^{k}$ in the positive semi-definite cone. Its closed-form solution can be simply realized by truncating all the information in the negative eigenvalue part. More precisely, for a symmetric matrix $X$ with eigen-decomposition $X=V \Lambda V^{\top}$, we define an eigenvalue hard thresholding (EVHT) operator as follows:
\begin{equation}\label{evht operator}
\mathcal{T}_E(X)=V \Lambda_+ V^{\top}, \ \ \  \Lambda_+=\text{diag}(\max\{\Lambda(i,i), 0\}_{i=1}^n).
\end{equation}
With this definition, the second step of \eqref{ReconModelAL_ alg} can be calculated as follows:
\begin{equation}\label{step2}
\hat{C}^{k+1}= \zeta\circ\mathcal{T}_E\circ\iota (\hat{B}^{k+1}+ H_2^{k}) .\\
\end{equation}

We summarize the above iterative method in Algorithm \ref{alg:ReconModel_B}. After finding $\hat{B}$ by Algorithm~\ref{alg:ReconModel_B}, eigen-decomposition of  $\iota(\hat{B})$ can ultimately reconstruct the coordinates of points $\{\bm{x}_1,\bm{x}_2,\ldots \bm{x}_\ell\}$ from the incomplete distance $\DO$.

\begin{algorithm}

\caption{Augumented Lagrangian method  to solve \eqref{eqn:ReconModel_B3}}
\label{alg:ReconModel_B}
\begin{algorithmic}

\STATE{\textbf{Initialization.}} Set $C^0=0,H_1^0=0,H_2^0=0,E(0)=2e+10,E(1)=1e+10,\epsilon=1e-6$. Set $\mu_1=10,\mu_2= 5$.

\WHILE{ ( $|E(k)-E(k-1)| /E(k) \geq \epsilon$) }

\STATE{\textbf{1.}} Solve  $ (\mu_1 \mathcal{A}^*\mathcal{A}+\mu_2) \hat{B}^{k+1}= \mu_1 \mathcal{A}^*(\TDO-H_1^{k}) +\mu_2(\hat{C}^k-H_2^{k})- \zeta(I_n)$.

\STATE{\textbf{2.}}  $\hat{C}^{k+1}= \zeta\circ\mathcal{T}_E\circ\iota (\hat{B}^{k+1}+ H_2^{k})$.

\STATE{\textbf{3.}}  $H_1^{k+1}=H_1^{k}+  \mathcal{A}\hat{B}^{k+1}-\TDO$.

\STATE{\textbf{4.}}  $H_2^{k+1}=H_2^{k}+   \hat{B}^{k+1} - \hat{C}^{k+1}$.

\STATE{\textbf{5.}}  Let $\displaystyle E(k+1) = \Tr(\iota(\hat{B}^{k+1})) + \frac{\mu_1}{2}\| \mathcal{A}\hat{B}^{k+1}-\TDO \|_2^2+ \frac{\mu_2}{2}\| \hat{B}^{k+1} - \hat{C}^{k+1} \|_2^2$.
\ENDWHILE
\end{algorithmic}
\end{algorithm}

\begin{remark}
In Algorithm \ref{alg:ReconModel_B}, Step 2 is usually time consuming as the full eigen-decomposition has time complexity $\O(\ell^3)$. Because ultimately the solution of $B$ and $C$ should be low rank, it is safe to only carry out the largest $m$ eigenfunctions by setting $m$ is relatively big enough comparing with the rank of $B$ but also small in terms of the matrix size. For instance, in our later numerical experiments for 3D point clouds, we use the MATLAB function "eigs" to find the largest $20$ eigenvalues and truncate all the other parts. In this way, the time consumption of Step 2 can be apparently decreased, and numerical results suggest no difference between this simplification and calculating full eigen-decomposition.
\end{remark}
\begin{remark}
\label{rk:LowrankD}
In addition, we would like to remark that the distance matrix is also a low rank matrix provided by $\text{Rank}(D) \leq \text{Rank}(B) +2 $ due to the relationship between $B$ and $D$ satisfying $D =  - 2B + \text{diag}(B) \bm{1}^{\top}  + \bm{1}  \text{diag}(B)^{\top}$, where $\text{diag}:\RR^{\ell\times \ell} \rightarrow \RR^{\ell \times 1}$ is the linear operator to project diagonal component of $B$ to a column vector. Therefore, we can also directly reconstruct $D$ based on its available information. This leads to a low rank minimization model as follows:
\begin{equation}\label{eqn:ReconModel_D}
\min\limits_{\tilde{D}\in\mathbb{S}} \|\tilde{D}\|_*  \quad  \text{s.t.}  \quad   \tilde{D}_{ij} = D_{ij}, \quad (i,j)\in\Omega  \\
\end{equation}
The above model is exactly the same as the general matrix completion theory discussed in \cite{CandesRe2008}, which also indicates that the minimum desired information of $D$ for successful reconstruction is proportion to the matrix rank. Note that the rank of $D$ is higher than the rank of matrix $B$. Plus, the less flexibility of $B$ due to the positive semi-definite condition can also help to reduce the requirement of valid distance information for successful reconstruction. Therefore, we would expect that model \eqref{eqn:ReconModel_D} requires a higher rate of valid distance than the proposed model \eqref{eqn:ReconModel_B2}. In fact, our numerical experiments in Section~\ref{sec:experiments} also verify this observation.
\end{remark}

\section{Solving PDEs on Manifolds from Incomplete Distance}
\label{sec:SolvePDEs}
In this section, after a brief review of the moving least square (MLS) method~\cite{Liang:CVPR2012,liang2013solving} and the local mesh (LM) method~\cite{lai2013local} for solving PDEs on manifolds represented as point clouds, we demonstrate our proposed methods of solving an elliptic eigenvalue problem and a nonlinear hyperbolic equation on manifolds represented as incomplete inter-point distance information.

\subsection{Approximating differential operators from reconstructed local coordinates}
\label{subsec:ReviewMLSLM}
For traditional implicit or triangle mesh representation of surfaces in $\mathbb{R}^3$, implicit methods, parameterization methods, finite difference methods and finite element methods~\cite{Osher:88,pinkall1993computing,Bertalmio:2000,taubin2000geometric,Bertalmio:2002,Meyer:2003,stam2003flows,spira2007geometric,wang2007brain,lui2008variational,Brandman:2008JSC,xu2004convergent,Reuter:06,lai2011framework,dziuk2013finite} have been proposed to solve PDEs on surfaces. However, it is not straightforward to use these methods in our case as the incomplete distance data structure has no global mesh structure and the ground manifold could be with dimension higher than two and in a high co-dimensional ambient space.
As the definitions of differential operators on manifolds are coordinate invariant, we can approximate differential operators point-wisely based on its local coordinate reconstruction obtained from the matrix completion method proposed in Section~\ref{sec:CoordinateRecon}. Namely, each given point, its reconstructed local neighborhood can be viewed as a point cloud, therefore, we apply moving least square (MSL) method~\cite{Liang:CVPR2012,liang2013solving} or local mesh (LM) method~\cite{lai2013local} to approximate the differential operators. These methods can achieve high order accuracy and enjoy more flexibility. They can be applied to manifolds with arbitrary dimensions and codimensions. To make the paper self-sufficient, we briefly discuss these two methods in this section.

Denote the index set of the KNN of the $i$-th point as $N(i)$ and write $\mathcal{X}(i) = \{\bm{x}_k\in\RR^p~|~ k\in N(i)\}$ as the reconstructed Euclidean  coordinates of the neighborhood $N(i)$. Based on $\mathcal{X}(i)$, its tangent space and normal space can be determined by standard principle component analysis (PCA)~\cite{jolliffe2002principal}, which is provided by eigensystem of the  covariance matrix $P_i$ of $N(i)$ defined on:
\begin{equation}\label{covmatrix}
P_i=\sum\limits_{k \in N(i)} (\bm{x}_k-c_i)^T (\bm{x}_k-c_i).
\end{equation}
where $c_i=\frac{1}{|N(i)|}\sum\limits_{k \in N(i)}\bm{x}_k$ is the centroid of $\mathcal{X}(i)$. If the intrinsic dimension of the manifold is $d$, then the jump of $P_i$'s eigenvalues guides the splitting $\RR^p = \T_i \oplus \N_i$. Here $\T_i$ represents the tangent space spanned by $\{e_i^1,\cdots,e^d_i\}$ corresponding to the $d$ largest eigenvectors of $P_i$, and $\N_i$ represents the normal space spanned by the rest of the eigenvectors of $P_i$.  As $\{e_i^1,\cdots,e_i^p\}$ forms an orthonormal basis near the point $\bm{x}_i$, a new {\it virtual coordinates} of $\bm{x}_k$ can be obtained by $(\bm{u}_k, \bm{v}_k) = (\{\langle \bm{x}_k, e_i^\alpha\rangle\}_{\alpha=1}^d,\{\langle \bm{x}_k, e_i^\beta\rangle\}_{\beta={d+1}}^p)$. Therefore, the manifold structure near $\bm{x}_i$ can be approximated by a degree 2 polynomial map $\displaystyle \bm{Z}_i:\RR^d\rightarrow \RR^{p-d}, \bm{Z}_i(\bm{u}) = (\sum_{|\bm{\xi}|\leq 2}\bm{c}^1_{\bm{\xi}} \bm{u}^{\bm{\xi}}, \cdots,\sum_{|\bm{\xi}|\leq 2}\bm{c}^{p-d}_{\bm{\xi}} \bm{u}^{\bm{\xi}})$ whose coefficients are given by the following moving least square problem:
\begin{equation}
\label{eqn:MLS_PC}
\min_{\bm{c}^1_{\bm{\xi}},\cdots,\bm{c}^{p-d}_{\bm{\xi}}}\sum\limits_{k \in N(i)} w(\|\bm{x}_k - \bm{x}_i\|) \|\bm{Z}_i(\bm{u}_k) - \bm{v}_k\|^2,
\end{equation}
where the standard multi-index notation is used for $\bm{u}^{\bm{\xi}} = u_1^{\xi_1}\cdots u_d^{\xi_d}$ and $w(\cdot)$ is a weight function whose typical choice can be $w(d)=\exp(-\frac{d^2}{h^2})$ and $\displaystyle h=\max_{k \in N(i)}\|\bm{x}_k-\bm{x}_i\|$. Similarly, any function $f$ or vector filed $\bm{V}$ defined on $\mathcal{X}(i)$ can also be interpolated as a degree 2 polynomial $f(\bm{u})$ or $\bm{V}(\bm{u})$. Overall, we obtain a local polynomial interpolation of $\mathcal{X}(i)$ as $\bm{Z}_i(\bm{u})$ and a local polynomial interpolation of function defined on $\mathcal{X}(i)$ as $f(\bm{u})$. Therefore, the metric tensor $G = (g_{st})$ near $\bm{x}_i$ can be obtained by $g_{st} = \delta_{st} + \langle \frac{\partial{\bm{Z}_i}}{\partial u_s}, \frac{\partial{\bm{Z}_i}}{\partial u_t} \rangle$. This enables us to approximate intrinsic differential operators at $\bm{x}_i$ such as:
\begin{equation*}
\nabla_\M f =g^{st} \frac{\partial f}{\partial u_t} \partial_{u_s}, \quad
\mathrm{div}_\M \bm{V} = \frac{1}{\sqrt{g}}  \frac{\partial}{\partial u_s} \left( \sqrt{g}~V^s\right),\quad
\Delta_\M f= \frac{1}{\sqrt{g}}  \frac{\partial}{\partial u_s} \left( \sqrt{g} g^{st} \frac{\partial f}{\partial u_t} \right).
\end{equation*}
where  $g = det(G)$, $(g^{st}) = G^{-1}$ and the Einstein summation is used. Other differential operators can also be approximated using the similar method. After that, a finite difference type of method can be applied to solve differential equations on the given data.
We refer~\cite{Liang:CVPR2012,liang2013solving} for more details about this approach.

The local mesh method proposed in \cite{lai2013local} is another way to approximate the above differential operators from the virtual local coordinates. Namely, at $i$-th point, a local connectivity can be constructed for the projection image of $\mathcal{X}(i)$ on $\T_i$ through the standard Delaunay triangulation. Thus, this connectivity can be directly inherited on $\mathcal{X}(i)$. If we write the simplex of the first ring of $\bm{x}_i$ by $\mathcal{R}(i)$, the differential operators at each point can be approximated by weighted average as follows:
\begin{equation*}
\displaystyle \nabla_{M}f(\bm{x}_i)=\frac{\sum_{S\in\R(i)}|S|\nabla_{S}f(\bm{x}_i)}{\sum_{S\in\R(i)}|S|},\quad
\mathrm{div}_{M}\overrightarrow{V}(p_i)=\frac{\sum_{S\in\R(i)}|S|\mathrm{div}_{S}\overrightarrow{V}(\bm{x}_i)}{\sum_{S\in\R(i)}|S|}
\end{equation*}
Moreover, a finite element type of method can be applied to estimate the mass matrix and stiffness matrix. Therefore, differential equations like Laplace-Beltrami eigenvalues problems can be solved. We refer~\cite{lai2013local} for more detailed discussion about construction of mass matrix and stiffness matrix.

We remark that the above constructions are conducted point-wisely. As long as a differential operator is well-defined on a manifold, namely, its definition does not depend on the choice of local coordinate, then the above procedure can consistently approximate the desired differential operators. A different strategy that can also be considered is to reconstruct coordinates first using the proposed Gram matrix completion algorithm, then approximate differential operators based on the global  reconstructed coordinates. We would like to point out that the proposed methods of approximating differential operators based on local reconstruction enjoys advantages of computation efficiency and memory consumption. In fact, the most time-consuming part of the coordinate reconstruction Algorithm~\ref{alg:ReconModel_B} is the eigenvalue hard thresholding. Consider a incomplete distance data with $n$ nodes, the complexity of each eigenvalue thresholding step is $\O(n^2m)$ for global reconstruction if only the largest $m$ eigenvalues are computed. In our local reconstruction strategy, we only reconstruct coordinates of $\ell$ nearest points near $\bm{x}_i$ with complexity
$\O(n\ell^2m)$. When $n$ is very large (such as $16002$ used in our experiments) and $\ell$ is small ( $\leq$ 30 in our experiments ), the local reconstruction strategy is apparently much more efficient than global reconstruction strategy. Our numerical experiments discussed in Section~\ref{sec:experiments} also support this point (see Table~\ref{tab:TimeLocalVSGlobal}). Moreover, the local reconstruction method supports parallel computation which can further reduce time consumption. Additionally, the memory consumption of global reconstruction is $\O(n^2)$ and the matrix is not sparse, which often crashes the program due to insufficient memory for large size data. However, our local reconstruction strategy only requires the memory $\O(\ell^2)$ which avoids the problem of exceeding memory limit.

\subsection{Solving differential equations based on incomplete distance}
\label{subsec:solvePDE}
The first PDE we consider is an elliptic eigenvalue problem of the Laplace-Beltrami (LB) operator $\Delta_\M$.
The LB operator is self-adjoint and elliptic, so its
spectrum is discrete. We denote the eigenvalues of $-\Delta_{\M}$
as $0=\lambda_0<\lambda_1<\lambda_2<\cdots$ and the corresponding
eigenfunctions as $\phi_0, \phi_1,\phi_2,\cdots$ satisfying the following equations~\cite{Chavel:1984}:
\begin{eqnarray}
\Delta_{\M} \phi_k=-\lambda_k\phi_k,  \quad k=0,1,2,\cdots.
\label{eqn:lb_closesurf}
\end{eqnarray}
The set of LB eigenfunctions $\{\phi_k\}_{k=0}^\infty$ forms an orthonormal basis of the space of $L^2$ functions on $\M$. The set
$\{\lambda_k,\phi_k\}_{k=0}^\infty$ is called LB eigensystem of $\M$.
Due to the intrinsic definition of the LB operator $\Delta_{\M}$, the induced LB eigensystem
$\{\lambda_k, \phi_k\}_{k=0}^ \infty$ is also completely intrinsic to the manifold geometry and provides an intrinsic and systematic characterization of the underlying manifold geometry~\cite{Berard:1994}.
Recently, there have been increasing interests in using the LB eigen-geometry for 3D shape analysis as well as point clouds analysis~\cite{Reuter:06,Levy:2006IEEECSMA,Vallet:2008CGF,Shi:08a,Peter:08,Sun:2009SGP,Lai:2010CVPR,Bronstein:2010CVPR,Liang:CVPR2012,lai2017multi}. Therefore, it would be also important to consider solve the LB eigensystem for manifolds represented as incomplete inter-point distance, then all existing methods of conducting data analysis using the LB eigensystem can also be adapted to the incomplete distance data structure. Our numerical solvers are based on the numerical solvers for on manifolds represented as point clouds. First, we apply the matrix completion model for point-wise local coordinate reconstruction. Either MLS or LM can be then applied to obtain the discretized matrix form of the equation, where the only step we need to conduct is to use the reconstructed local coordinates for the LB operator discretization as we discussed in the section~\ref{subsec:ReviewMLSLM}. As the LB operator is invariant under different choices of the local coordinates, the rigid motion ambiguity from the local coordinate reconstruction will not introduce inconsistency across different coordinate patches. This guarantees that our method can provide satisfactory numerical results as we illustrate in the numerical section. We would like to remark that it is relatively straightforward to have the local consistency of the differential operator as long as the available distance information satisfying certain incoherence condition as we discussed in our theoretical validation work~\cite{lai2017Exact}. However, the stability of the discretization is still open and will be explored in our future work.

The second equation we consider is the Eikonal equation, a special type of nonlinear hyperbolic PDE on manifolds. This equation is used to characterize the geodesic distance, an intrinsic measurement between two points, on a manifold.
The Eikonal equation for the distance map $d$ to a given set $\Gamma$ on $\M$ can be stated as follows:
\begin{equation}\label{eqn:Eikonal}
\begin{cases}
| \nabla_\M d(x)| = 1 \\
d(x) =0,\ x \in \Gamma \subset \M
\end{cases}
\end{equation}
when $\Gamma$ only includes a single point $p$, the distance map can be denoted as $d_p$. In practice, provided points coordinates and mesh structure, model \eqref{eqn:Eikonal} can be solved by fast marching and fast sweeping methods \cite{sethian1996fast,kimmel1998computing,zhao2005fast}. For manifolds represented as point clouds, a local mesh method based on fast marching is proposed in~\cite{lai2013local}. As long as the local coordinate reconstruction is obtained using the proposed matrix completion model, then a local Delaunay mesh structure can be constructed. After that, we repeat the local mesh method to conduct the fast marching as this approach only depends on the distance information of the first ring structure of the local mesh, which is again rigid motion invariant and will not be affected due to the rigid motion ambiguity from the local coordinate reconstruction.

\section{Applications on Manifold Reconstruction and Dimension Reduction}
\label{sec:MfdStitching}
As we mentioned before, solutions of differential equations on manifolds can provide global information for understanding data structure. Thus, PDEs can be viewed as ``bridges" linking between local information and global information. This shares the same sprit as ``think globally, fit locally" discussed in~\cite{saul2003think} although no PDEs are involved there. Therefore, without global reconstruction, some global analysis of point clouds such as pattern recognition, comparisons and classification can be further considered for data represented as incomplete distance. In this section, we illustrate results on applications of manifold reconstruction and dimension reduction based on solutions of the Laplace-Beltrami eigenvalue problem on incomplete distance data.
\subsection{Manifold reconstruction from distance via patch stitching using LB eigenfunctions}
The patch stitching problem is to reconstruct the coordinate of a point cloud $P \in \mathbb{R}^{n \times d}$ from coordinates of $L$ subsets (referred as patches) $\{\Omega_j\}_{j=1}^{L}$, where $\Omega_j$ denotes the index set for patch $j$. In practice, for each patch $\Omega_j$, the local coordinates $Q_j\in\RR^{n_j\times d}$, which are provided as input information or are reconstructed from distance information, have the same pair-wise distance and geometry as the restricted global coordinate $P_{\Omega_j}$ with possible rigid motion difference. In other words, $P_{\Omega_j}= Q_j R_j + \bm{1}_j b_j$, where $R_j \in O(d) = \{ R\in\RR^{d \times d} ~|~ R^{\top} R = I_d \}$ is an orthogonal matrix for rotation and reflection, $\bm{1}_j \in \mathbb{R}^{n_j \times 1}$ is a vector with all element equal to 1, and $b_j \in \mathbb{R}^{1\times d}$ is the translation vector. A straightforward approach to solve $\{R_j\}$ and $\{b_j\}$ is to minimize the quadratic loss:
\begin{equation}
\label{eqn:qualoss}
\min\limits_{P,\{R_j \in O(d), ~b_j \in\RR^{d}\}_{j=1}^{L}} \sum_{j=1}^{L} \|P_{\Omega_j}- Q_j R_j - \bm{1}_j b_j\|_2^2,
\end{equation}
Such model is a non-convex problem since $\{R_j\}$ are required to optimize over the non-convex domain of orthogonal transforms. In \cite{cucuringu2012sensor}, the authors proposed a three-stage patch stitching method from the local coordinate by synchronizing the reflection, rotation and translation successively. However, the quality of synchronization is strongly relies on the proportion of the overlapping index. More recently, a spectral relaxation method is proposed to solve the above problem in \cite{chaudhury2015global}. This method essentially relax the set of orthogonality constraints and perform well if neighborhood patches has reasonable enough overlapping points.

Here, we propose an alternative method to reconstruct the global coordinate from local coordinate, or merely incomplete distance information of each patch. Based on solution of  LB eigen-problem incomplete distance information discussed in Section \ref{subsec:solvePDE},  we consider the LB eigen-system as a "bridge" to connect the information of different patches. This idea of this new method is essentially align each local patch in the spectral domain instead of the original Euclidean coordinate domain. As the global information has been captured from the LB eigen-system, the overlapping of different patches is not directly required here. This approach is inspired by approximating coordinates using LB eigenfunctions. For a given point cloud $P \in \mathbb{R}^{n \times d}$, it can be approximated by the linear combination of the LB eigenfunctions as $P \approx \Phi \bm{\alpha}$. Here, $\Phi = (\phi_1,\cdots,\phi_N) \in \mathbb{R}^{n \times N} $ represents the first $N$ LB eigenfunctions $\{\phi_i\}_{i = 1}^N$ and $\bm{\alpha} = (\alpha_1, \cdots, \alpha_d) \in\RR^{N\times d}$ represents coefficients of coordinates as column vectors. Since the mass matrix is generally unknown without the prior knowledge of the structure information, the coefficients matrix $\bm{\alpha}$ can be obtained in the following least square sense:
\begin{equation}\label{solvealpha1}
\min_{\bm{\alpha}}  \|P -  \sum\limits_{i=1}^{N} \phi_i \alpha_i \|_F^2 = \min_{\bm{\alpha}}  \|P -  \Phi \bm{\alpha} \|_F^2,
\end{equation}
Consider the relationship between restricted global coordinate $\{P_{\Omega_j}\}$ and the input local coordinate $\{Q_j\}$, we already have that for any $1 \leq j \leq L$, $P_{\Omega_j}= Q_j R_j + \bm{1} b_j$. Therefore the norm $ \|P -  \Phi \bm{\alpha} \|_2^2$ in local patch $Q_j$ can be equivalently written as $ \|P_{\Omega_j} -  \Phi_{\Omega_j} \bm{\alpha} \|_2^2= \|Q_j R_j + \bm{1}_j b_j -  \Phi_{\Omega_j} \bm{\alpha} \|_2^2= \|Q_j  + \bm{1}_j b_j R_j^{\top}-  \Phi_{\Omega_j} \bm{\alpha} R_j^{\top}\|_2^2$. Since $b_j R_j^{\top}$ is still a translation vector in $\mathbb{R}^d$, the norm can be simplified as $\|Q_j  -  \Phi_{\Omega_j} \bm{\alpha} R_j - \bm{1}_jb_j \|_2^2$ without loss of generality. As a result, with the given $\Phi$ and all the local coordinate $\{Q_j\}$, we propose the following model for stitching point clouds $\{Q_j\}_{j=1}^{L}$:
\begin{equation}\label{solvealpha2}
\min_{\bm{\alpha}\in\RR^{N\times d},\{R_j \in O(d), ~b_j \in\RR^{d}\}_{j=1}^{L}} \E(\bm{\alpha}, R_j, b_j) =  \frac{1}{2} \sum\limits_{j=1}^{L}   \|Q_j -  \Phi_{\Omega_j} \bm{\alpha} R_j- \bm{1}_jb_j \|_F^2,
\end{equation}
Once $\{R_j\}$ and $\{b_j\}$ are obtained, we can find all the adjusted coordinates $P_{\Omega_j}= Q_j R_j + \bm{1}_j b_j$, which leads to the global coordinate $P$ ultimately.

Different from the synchronization method and the spectral relaxation method considered in \cite{cucuringu2012sensor,chaudhury2015global}, we design the following method to solve the nonconvex problem  \eqref{solvealpha2} by updating  $\bm{\alpha}$, $\{R_j\}$ and $\{b_j\}$ iteratively:
\begin{equation}\label{gradient descent 1}
\begin{cases}
\bm{\alpha}^{k+1}= \bm{\alpha}^k - (\sum \limits_{j=1}^{L} \Phi_{\Omega_j}^{\top}\Phi_{\Omega_j})^{-1} \nabla _{\bm{\alpha}}\E(\bm{\alpha}^k, R_j^k, b_j^k) ,\\
R_j^{k+1}= \mathcal{K}(R_j^{k}, \nabla _{R_j}\E(\bm{\alpha}^{k+1}, R_j^k, b_j^k)), \qquad \forall 1 \leq j \leq L,\\
\displaystyle b_j^{k+1} = b_j^{k} -  \frac{1}{\text{Card}(\Omega_j)} \nabla _{b_j}\E(\bm{\alpha}^{k+1}, R_j^{k+1}, b_j^k) , ~\quad \forall 1 \leq j \leq L.\\
\end{cases}
\end{equation}
where we update $\bm{\alpha}$ and $\{b_j\}$ using Newton's method and it is straightforward to check
\begin{equation}
\label{eqn:gradientofstitching}
\begin{cases}
\nabla _{\bm{\alpha}}\E(\bm{\alpha}^k, R_j^k, b_j^k) = - \sum \limits_{j=1}^{L}  \Phi_{\Omega_j}^{\top}(Q_{j} -  \Phi_{\Omega_j} \bm{\alpha}^k R_j^{k}- \bm{1}_jb_j^{k} ) (R_j^{k})^{\top}, \\
\nabla _{R_j}\E(\bm{\alpha}^{k+1}, R_j^k, b_j^k) = - (\bm{\alpha}^k)^{\top} \Phi_{\Omega_j}^{\top}(Q_{j} -  \Phi_{\Omega_j} \bm{\alpha}^k R_j^{k}- \bm{1}_j b_j^{k} ),\\
\nabla _{b_j}\E(\bm{\alpha}^{k+1}, R_j^{k+1}, b_j^k) = -\bm{1}_j^{\top}  (Q_{j} -  \Phi_{\Omega_j} \bm{\alpha}^k R_j^{k}- \bm{1}_j b_j^{k} ).\\
\end{cases}
\end{equation}

We next describe a method of updating the orthogonality constrained variable $\{R_j\}$. The operator $\mathcal{K}$ proposed in \cite{wen2013feasible} is designed for searching the gradient descent direction on the Stiefel manifold, the geometric description of the set of orthogonal matrices. To realize the operator $\mathcal{K}(\nabla R_j^{k})$, we first define the skew-symmetric operator as
\begin{equation}\label{skew-symmetric}
\begin{aligned}
G_j^k &= \nabla _{R_j}\E(\bm{\alpha}^{k+1}, R_j^k, b_j^k) &\\
S_j^k &= G_j^k (R_j^{k})^{\top} - R_j^{k}(G_j^k)^{\top}, & \forall 1 \leq j \leq L,\\
\end{aligned}
\end{equation}
\noindent then the new trial point $R_j^{k+1}$ satisfing $R_j^{k+1} (R_j^{k+1})^{\top} = I$ can be generated by
\begin{equation}\label{gradient descent so3}
R_j^{k+1} = \mathcal{K}(R_j^{k}, G^k_j) = (1+\frac{\delta}{2} S_j^k)^{-1} (1-\frac{\delta}{2} S_j^k)  R_j^{k},\\
\end{equation}
where the step size $\delta$ can be obtained by setting fixed value or line search methods. More detailed discussion about the orthogonality preserving property and convergence of this approach can be found in \cite{wen2013feasible}. With the above explanation of \eqref{gradient descent 1}, we summarize an algorithm of solving \eqref{solvealpha2} as Algorithm \ref{Algorithm stitching}.

\begin{algorithm}[h]

\caption{Gradient Descent method  to solve \eqref{solvealpha2}}
\label{Algorithm stitching}
\begin{algorithmic}

\STATE{\textbf{Initialization:}} Set the initial values such that $\bm{\alpha}^0=0, b_j^0 =0, R_1^0=I_{d \times d}$. By roughly estimating the orthogonal transform matrix $R_{i,j}$ between neighbourhood patches $\Omega_i$ and $\Omega_j$ with overlapping points, one can set the initial guess of $R_j^0 = R_i^0 R_{i,j}$. Using the broad first search (BFS) scheme, the rotation matrix of all connected patches can be estimated. If a patch $\Omega_i$ is isolated to any visited patches in the BFS algorithm, set $R_i^0=I_{d \times d}$.

\WHILE{ $\frac{|E(k)-E(k-1)|}{E(k)} >  \epsilon$}

\STATE{\textbf{1.}}  $\displaystyle \bm{\alpha}^{k+1}= \bm{\alpha}^k + (\sum \limits_{j=1}^{L} \Phi_{\Omega_j}^{\top}\Phi_{\Omega_j})^{-1} (\sum \limits_{j=1}^{L}  \Phi_{\Omega_j}^{\top}(Q_{j} -  \Phi_{\Omega_j} \bm{\alpha}^k R_j^{k}- \bm{1}_j b_j^{k} ) (R_j^{k})^{\top})$.  \vspace{0.1cm}

\STATE{\textbf{2.}}  $\displaystyle R_j^{k+1} = (1  + \frac{\delta}{2} S_j^k)^{-1}(R_j^{k}- \frac{\delta}{2} S_j^kR_j^{k})$ where $ G_j^k =  - (\bm{\alpha}^k)^{\top} \Phi_{\Omega_j}^{\top}(Q_{j} -  \Phi_{\Omega_j} \bm{\alpha}^k R_j^{k}- \bm{1}_j b_j^{k} )$, $S_j^k = G_j^k (R_j^{k})^{\top} - R_j^{k}(G_j^k)^{\top} ~\forall 1 \leq j \leq L $ and $\delta$ is the step size which can be solved by line search method.   \vspace{0.1cm}

\STATE{\textbf{3.}}  $\displaystyle b_j^{k+1} = b_j^{k} + \frac{1}{\text{Card} (\Omega_j)} \bm{1}_j^{\top}  (Q_{j} -  \Phi_{\Omega_j} \bm{\alpha}^k R_j^{k}- \bm{1}_j b_j^{k} ), ~\forall 1 \leq j \leq L$,

\STATE{\textbf{4.}}  Calculate $\displaystyle E(k+1)= \E(\bm{\alpha}^{k+1}, R_j^{k+1}, b_j^{k+1}) $ for the stopping criteria.

\ENDWHILE

\STATE{\textbf{5.}} Apply $P_{\Omega_j}= Q_j R_j + \bm{1}_j  b_j$ and combine all the $P_{\Omega_j}$ to find the $P$.
\end{algorithmic}
\end{algorithm}

Based on the above method, we therefore propose a global coordinate reconstruction algorithm from merely random missing distance in each local patches. Given the local patches $\{\Omega_j\},1 \leq j \leq L$ with patch size not larger than $K$, we assume the local distance information include $r \%$ of the local distance of each patches. The full reconstrution procedure can be summarized as follows:

\begin{description}
\item[Step 1]  Using the method discussed in Section~\ref{subsec:solvePDE} with local mesh reconstruction method to compute the first $N$ Laplace-Beltrami eigenfunctions $\Phi$.
\item[Step 2]   Using Algorithm \ref{alg:ReconModel_B} to generate the virtual local coordinate $\{Q_j\},1 \leq j \leq L$ for each patch.
\item[Step 3]   Using Algorithm \ref{Algorithm stitching} to find the global coordinate $P$ from $\Phi$ and $\{Q_j\}$.
\end{description}

%
%
%
%
%
%

\section{Numerical Results}
\label{sec:experiments}
In this section, numerical tests are presented to illustrate the proposed methods for solving PDEs on manifold represented as incomplete inter-point distance information and the application to manifold stitching. First, we demonstrate the effectiveness of the proposed matrix completion method for coordinate reconstruction. We also illustrate the phase transition curve of successful reconstruction. Second, we test our method of computing the LB eigen-problem and illustrate that the proposed method of solving PDEs through local coordinate reconstruction is much more efficient than solving PDEs from global coordinate reconstruction. In addition, we also test the problem of solving Eiknoal equations and illustrate the accuracy improvement of this approach to the canonical Dijkstra method. As applications, we also demonstrate preliminary results on point clouds reconstruction using global information from LB eigenfunctions and its potential application to dimension reduction problem. All numerical experiments are implemented by MATLAB in a PC with a 32G RAM and 16 dual-core 2.7 GHz CPUs.

\subsection{Coordinates reconstruction from incomplete distances}
The local coordinate reconstruction method based on matrix completion is a crucial part of solving differential equations on a manifold represented as incomplete inter-point distance data. To demonstrate the effectiveness of the method proposed in Section \ref{sec:CoordinateRecon}, we first test this method for global coordinate reconstruction, although our numerical method for solving differential equations only conduct coordinate reconstruction for a certain neighborhood of each point. Starting with a full distance matrix $D \in \mathbb{R}^{n \times n}$, the actual number of freedom in $D$ is $(n^2-n)/2$ as the diagonal elements of $D$ are $0$ and $D$ is symmetric. We randomly chose $\gamma (n^2-n)/2, (\gamma \in[0,~1])$ number of entries of $D$ as the input incomplete distance matrix and apply  Algorithm \ref{alg:ReconModel_B} to find a Gram matrix $B$. After that, global coordinates can be consequently generated from $B$ by its eigen-decomposition as we described in Section \ref{subsec:MDS}.

\begin{figure}[h]
\centering
\begin{tabular}{c@{\hspace{2pt}}c}
\includegraphics[width=0.49\linewidth]{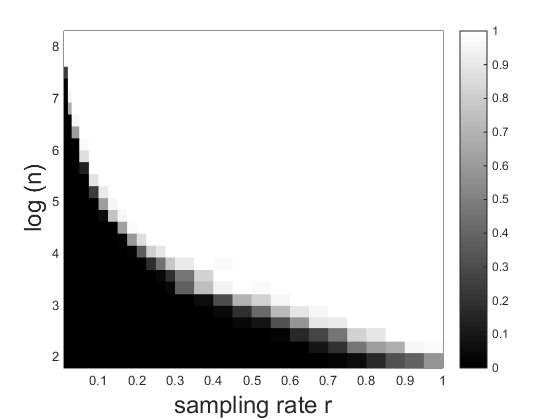}&
\includegraphics[width=0.49\linewidth]{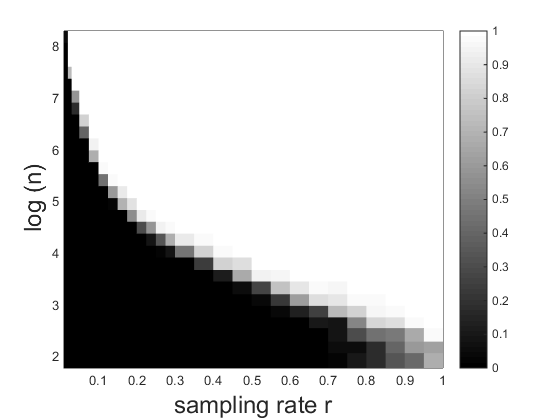}\\
\end{tabular}
\caption{Phase transition of successful reconstruction rate $\rho$ out of $50$ tests of reconstructing a uniformly sampled unit sphere. Left: successful reconstruction rate of the proposed model \eqref{eqn:ReconModel_B2} based on Gram matrix. Right: successful reconstruction rate based on distance completion model \eqref{eqn:ReconModel_D}.}
\label{fig:borderline1}
\end{figure}

Our first numerical test is conducted for checking successful reconstruction rate of reconstructing the input incomplete distance matrices. We test our model for different size matrices with fixed rank. By uniformly choosing $n$ points on the unit sphere, we obtained the corresponding pairwise Euclidean distance matrix $D_T$ and the associated Gram matrix $B_T$ as ground truth. After that, we randomly choose $\gamma n(n-1)/2$ entries of $D_T$ as available information and apply our reconstruction algorithm. We mean that the reconstruction is successful if the relative error between $B_T$ and the reconstructed Gram matrix under the Frobenius norm is less than $10^{-3}$. Given $n$ and $\gamma$, we run 50 tests and record the rate of successful reconstruction as $\displaystyle \rho = \frac{\# \mbox{ successful reconstruction}}{50}$. Based on the same input, Figure \ref{fig:borderline1} plots the phase transition of $\rho$ using Gram matrix completion model \eqref{eqn:ReconModel_B2} (the left image in Figure \ref{fig:borderline1})  and distance completion model \eqref{eqn:ReconModel_D} (the right image in Figure \ref{fig:borderline1}) in terms of the sampling rate $\gamma$ and the logarithm of total number of points $\log(n)$. Both images clearly show that for a fixed rank problem, successful reconstruction needs less portion of distance information as number of points increases. As we expect, the Gram matrix completion model \eqref{eqn:ReconModel_B2} has a larger domain of successful reconstruction than the model \eqref{eqn:ReconModel_D} using distance matrix completion. This is because the rank of Gram matrix is less than the rank of distance matrix and an additional semi-positive definite constraint is also impose for Gram matrix which further decrease the number of freedom. Therefore, we use model \eqref{eqn:ReconModel_B2} as our major local coordinate reconstruction tool to approximate differential operators in the rest of numerical experiments. 

\begin{table}[h]
\centering
{\footnotesize
\begin{tabular}{|c|ccccccc|}
\hline
\backslashbox{Data}{$\gamma$} &&  1\% &2\% & 3\% &5\%&10\% &20\% \\
\hline\hline
\multicolumn{1}{|c|}{\multirow{2}{*}{$S^2$} } &\multicolumn{1}{ |c|  }{$E_{B}$} & 7.157E-1 & 1.376E-3 & 4.791E-4 & 2.474E-4 &1.342E-5 & 4.262E-5\\
\cline{2-8}
& \multicolumn{1}{ |c|  }{$\rho$} & 0\% &92\%  &100\%  &100\%  &100\% &100\%\\
\hline\hline
\multicolumn{1}{ |c|  }{\multirow{2}{*}{Cow} }  & \multicolumn{1}{ |c|  }{$E_{B}$}  & 4.9427E-5 & 3.980E-4 & 1.837E-4 & 5.319E-5 & 1.4072E-5 & 2.155E-5\\
\cline{2-8}
& \multicolumn{1}{ |c|  }{$\rho$} & 100\% &100\%  &100\%  &100\%  &100\% &100\%\\
\hline\hline
\multicolumn{1}{ |c|  }{\multirow{2}{*}{Swiss roll} }  & \multicolumn{1}{ |c|  }{$E_{B}$} & 2.722E-4 &2.894E-4 &1.633E-4 &5.054E-5 &1.704E-5& 1.114E-5 \\
\cline{2-8}
& \multicolumn{1}{ |c|  }{$\rho$} & 100\% &100\%  &100\%  &100\%  &100\% &100\%\\
\hline
\end{tabular}
}
\caption{Rate of The successful reconstruction $\rho$ and the average relative error $E_{B}$ of the Gram matrix out of $50$ tests by the proposed model \eqref{eqn:ReconModel_B2} from distances with information availability rate $\gamma$. }
 \label{tab:tablerdip}
\end{table}

\begin{figure}[htp]
\centering
\begin{minipage}{0.24\linewidth}
\includegraphics[width=1\linewidth]{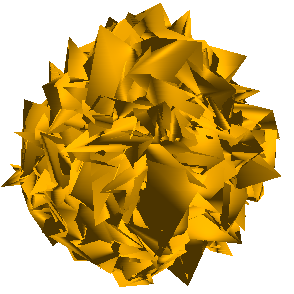}
\end{minipage}
\begin{minipage}{0.24\linewidth}
\includegraphics[width=1\linewidth]{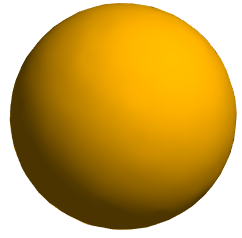}
\end{minipage}
\begin{minipage}{0.24\linewidth}
\includegraphics[width=1\linewidth]{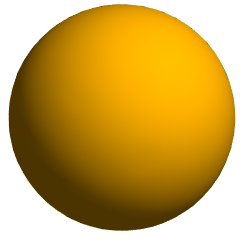}
\end{minipage}
\begin{minipage}{0.24\linewidth}
\includegraphics[width=1\linewidth]{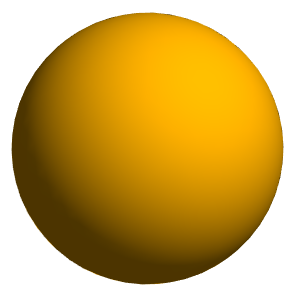}
\end{minipage}\hfill\\
\centering
\begin{minipage}{0.24\linewidth}
\includegraphics[width=1\linewidth]{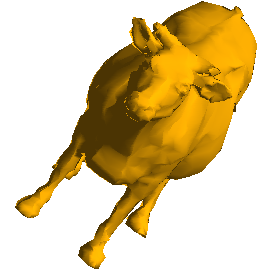}
\end{minipage}
\begin{minipage}{0.24\linewidth}
\includegraphics[width=1\linewidth]{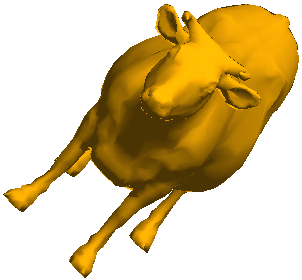}
\end{minipage}
\begin{minipage}{0.24\linewidth}
\includegraphics[width=1\linewidth]{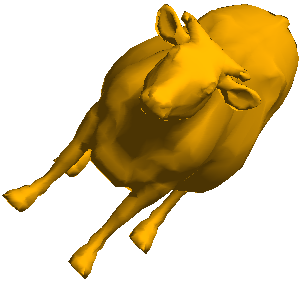}
\end{minipage}
\begin{minipage}{0.24\linewidth}
\includegraphics[width=1\linewidth]{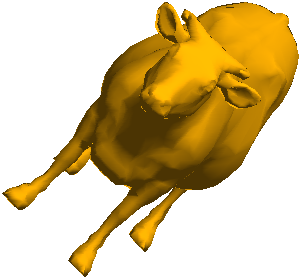}
\end{minipage}\hfill\\
\caption{Globally reconstructed coordinates of a unit sphere (1002 points) and a cow (2601 points) using proposed model \eqref{eqn:ReconModel_B2}. Images from left top to right bottom are reconstructed from distance with of 1\%, 2\%, 3\%, and 5\% available informatoin, respectively. }
\label{fig:randomball}
\end{figure}

We further test distance matrices from different data including point clouds sampled from the unit sphere(1002 points), a cow surface (2602 points), and a Swiss roll surface (2048 points) using our Gram matrix completion model. Table \ref{tab:tablerdip} reports the relative error of the Gram matrix and rates of successful reconstruction out of 50 tests for each $\gamma$. Figure \ref{fig:randomball} shows that with 1\% distance information, the reconstructed coordinates cannot formulate an outline of the unit sphere while the same rate of distance information can successfully reconstruct the cow shape. This is compatible with our observation showed in Figure \ref{fig:borderline1}, which indicates less percentage information is need for successfully reconstructing matrices with larger size. In addition, Table~\ref{tab:tablerdip} and Figure \ref{fig:randomball} also demonstrate both quantitatively and qualitatively that the proposed reconstruction model can provide very good coordinate reconstruction once the available distance information is sufficient.

\subsection{Solve LB eigenvalue problem from distance data with missing values}
\label{subsec:SovelLBeigs}

\begin{figure}[h]
\centering
\begin{minipage}{0.38\linewidth}
\centering{\footnotesize $\lambda = 20$}\\
\includegraphics[width=1\linewidth]{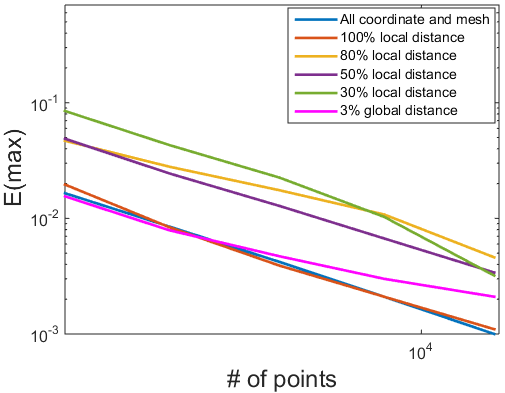}\\
\centering{\footnotesize $\lambda = 72$}\\
\includegraphics[width=1\linewidth]{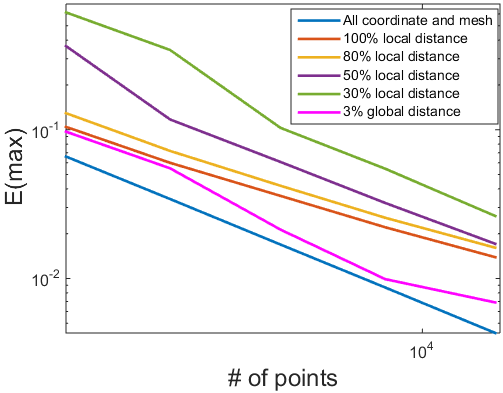}
\end{minipage}
\begin{minipage}{0.6\linewidth}
{\footnotesize
\begin{tabular}{|c|c|c|c|c|c|}
\hline
{\tiny $\#$ of points}& 1002&1962 & 4002 &7842&16002\\ \hline\hline
\multicolumn{6}{|c|}{finite element method with coordinates and mesh}\\ \hline
$\lambda=20$ & 0.0165 & 0.0085& 0.0042 & 0.0021 &0.0010\\ \hline
$\lambda=72$ & 0.0660 & 0.0342 & 0.0169 &0.0087 & 0.0043 \\
\hline\hline
\multicolumn{6}{|c|}{$\gamma = 100\%,\qquad  \ell=6$}\\ \hline
$\lambda=20$ & 0.0196 & 0.0084& 0.0039 & 0.0021 &0.0011\\ \hline
$\lambda=72$ & 0.1045 & 0.0599 & 0.0359 &0.0221 & 0.0139 \\
\hline\hline
\multicolumn{6}{|c|}{$\gamma = 80\%,\qquad  \ell=9$}\\ \hline
$\lambda=20$ & 0.0469 & 0.0280 & 0.0175 & 0.0108 &0.0046\\ \hline
$\lambda=72$ & 0.1292 & 0.0720 & 0.0420 &0.0256 & 0.0161 \\
\hline\hline
\multicolumn{6}{|c|}{$\gamma = 50\%,\qquad  \ell=18$}\\ \hline
$\lambda=20$ & 0.0488 & 0.0245 & 0.0128 & 0.0067 &0.0034\\ \hline
$\lambda=72$ & 0.3649 & 0.1172 & 0.0607 &0.0322 & 0.0171 \\
\hline\hline
\multicolumn{6}{|c|}{$\gamma = 30\%,\qquad  \ell=30$}\\ \hline
$\lambda=20$ & 0.0848 & 0.0431 & 0.0225 & 0.0103 &0.0032\\ \hline
$\lambda=72$ & 0.6146 & 0.3435 & 0.1029 &0.0547 & 0.0263 \\
\hline\hline
\multicolumn{6}{|c|}{3\% of global distance information ($\ell=6$ for MLS)}\\ \hline
$\lambda=20$ & 0.0155 & 0.0079& 0.0047 & 0.0030 &0.0021\\ \hline 
$\lambda=72$ & 0.0970 & 0.0550 & 0.0213 &0.0099 & 0.0069\\
\hline
\end{tabular}
}
\end{minipage}
\caption{Computation errors for LB eigenvalues $\lambda=20, 72$ on the unit sphere based on incomplete distance $D_{\gamma,\ell}$ with KNN information. Left: convergence curves. Right: Relative errors for input distance matrices with different sample size.}
\label{fig:RelErrorsforLBeigenvaluesOnS2}
\end{figure}

Our second numerical experiment illustrates the proposed method of solving the LB eigenvalue problem on manifolds represented as incomplete inter-point distance information. Considering a point cloud $\{p_1,\cdots, p_n\}$ uniformly sampled on a given manifold $\mathcal{M} \subset \RR^p$, we compute the associated distance matrix $D$ provided by pairwise Euclidean distance for this point cloud. Our numerical experiments illustrate the proposed method for computing the LB eigenvalue problem only based on the incomplete distance. More precisely, we first uniformly randomly choose a subset $\Omega_\gamma \subset \{(i,j) ~|~  1\leq j< i \leq n \}$, such that $|\Omega_\gamma| = \left\lceil \gamma n(n-1)/2 \right\rceil$.
Then incompletion distance information is provided as
\begin{equation*}D_{\gamma,\ell}  = \Big\{D(i, j)~|~ j\in N_{\ell}(i) ~\&~ (i,j) \in\Omega_\gamma, \quad i=1,\cdots,n \Big\}
\end{equation*}
where $N_{\ell}(i)$ denotes the index set of the nearest $\ell$ points to $p_i$. Thus, $D_{\gamma,\ell}$ roughly contains $\gamma \ell/n$ portion of $D$, which contains a very small portion of $D$ if $\ell\ll n$ as we showed in our experiments. In addition, we measure computation accuracy using the normalized error $\displaystyle E_{\max,k}= \max_{i}\frac{|\tilde{\lambda}_{k,i}-\lambda_k|}{|\lambda_k|}$, where $\tilde{\lambda}_{k,i}$ are numerical approximation from our method, and $i$ represents the multiplicity.

We first demonstrate our method on different data sampled from the unit sphere. It is known that the $k$-th LB eigenvalue of the unit sphere is given by $\lambda_k=k(k+1)$ with multiplicity $2k+1$. This allowss us to illustrate numerical accuracy of our methods on the unit sphere case. In the first test, we assume the $K$-nearest neighborhood (KNN) is provided and sample distance matrices for point clouds on the unit sphere by setting different values of $\gamma = 100\%, 80\%, 50\%, 30\%$  and with different KNN size $\ell = 6, 9, 18, 30$, respectively. 
Given a KNN size $\ell$, the settings of local sampling rate should above the border line as shown in Figure \ref{fig:borderline1} to guarantee the exact local coordinate reconstruction. Conversely, given a fixed sampling rate of local distance, we need to set $\ell$ large enough to guarantee the successful reconstruction.
With these input data, we apply the proposed strategy for approximating the LB operator only from local coordinates reconstruction and compute the first 100 LB eigenvalues. We typically chose approximation results of LB eigenvalues $\lambda=20$ and $\lambda = 72$ to conduct accuracy analysis by illustrating the numerical error $E_{\max,n}$ in the right table of Figure \ref{fig:RelErrorsforLBeigenvaluesOnS2}. In addition, we list results from global coordinate reconstruction based on $3\%$ of distance information. As a reference, we also list computation results from classical finite element method from all the exact coordinates of points and the mesh information, which can be experimentally regarded as the best computation result based on the whole distance information. It is clear that our method provides satisfactory approximation of the LB eigenvalue problem on the unit sphere. Convergence of our method is illustrated by the error curves showed in the left two plots in Figure \ref{fig:RelErrorsforLBeigenvaluesOnS2}.  Moreover, we also observe the second order convergence numerically as $E_{\max,k}$ is approximately reduced by half if the number of total point doubled. This is compatible with convergence behavior of the moving least square method for solving the LB eigenvalue problem on point clouds discussed in~\cite{Liang:CVPR2012,liang2013solving}.

\begin{table}[h]
\centering
{\footnotesize
\begin{tabular}{*{5}{|p{1.2cm}}|}
\hline
\multicolumn{5}{|c|}{number of points}\\ \hline
1002&1962 & 4002 &7842&16002\\ \hline
\multicolumn{5}{|c|}{$\gamma = 100\%,\qquad  \ell=6$}\\ \hline
0.26 & 0.51 & 1.01 & 2.03 &4.05\\ \hline

\multicolumn{5}{|c|}{$\gamma = 80\%,\qquad  \ell=9$}\\ \hline
2.28 & 5.60 & 11.17 & 22.28 &45.02\\ \hline

\multicolumn{5}{|c|}{$\gamma = 50\%,\qquad  \ell=18$}\\ \hline
4.03 & 8.09 & 16.14 & 32.44 &64.71\\ \hline

\multicolumn{5}{|c|}{$\gamma = 30\%,\qquad  \ell=30$}\\ \hline
15.13 & 30.19 & 60.42 & 120.95 &241.63\\ \hline \hline

\multicolumn{5}{|c|}{\textbf{global} reconstruction using 3\%  distance  ($\ell=6$ for MLS)}\\ \hline
2.09 & 9.86 & 40.13 & 154.40 & 597.06\\ \hline 
\end{tabular}
}
\caption{Comparisons of time consumption (minutes) of solving the LB eigenvalue problem based on local/global reconstruction methods.}
\setlength{\abovecaptionskip}{10pt}
 \label{tab:TimeLocalVSGlobal}
\end{table}
We also report total time consumption in Table~\ref{tab:TimeLocalVSGlobal} for computing the first 100 LB eigenvalues with different setting of $\gamma$ and $\ell$ in this test. Note that the major time-consuming part is local coordinates reconstruction based on incomplete distance information, whose complexity is dependent on the local size parameter $\ell$ and the total number of points $n$. Theoretically, our local coordinate reconstruction method generally has the time complexity $\mathcal{O}(\min(\ell^2nm,\ell^3 n))$, while the global coordinate reconstruction has time complexity $\mathcal{O}(n^2 m)$. Here, $m$ is denoted as the maximum number of eigenvalues used in the eigenvalue hard thresholding step of algorithm~\ref{alg:ReconModel_B}, which is typically chosen as 20 in our experiments. For large data sets satisfying $\ell \ll n$, local reconstruction method can save a huge amount of time. This advantage is also illustrated  in Table~\ref{tab:TimeLocalVSGlobal}, where local reconstruction method for solving LB eigenvalue problem is about 150 times faster than global reconstruction method if we choose $n = 16002$ and $\ell = 6$.

\begin{figure}[htp]
\centering
\begin{minipage}{0.38\linewidth}
\centering{\footnotesize $\lambda = 20$}\\
\includegraphics[width=1\linewidth]{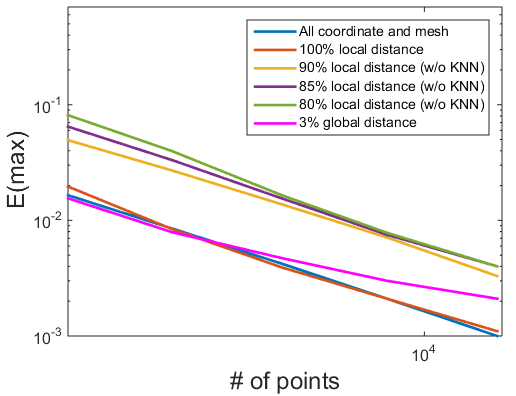}\\
\centering {\footnotesize $\lambda = 72$}\\
\includegraphics[width=1\linewidth]{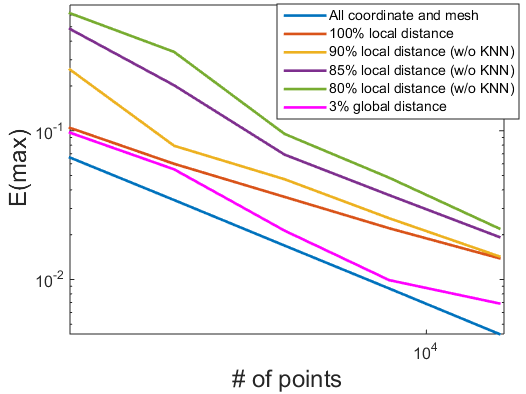}
\end{minipage}
\begin{minipage}{0.61\linewidth}
{\footnotesize
\begin{tabular}{|c|c|c|c|c|c|}
\hline
{\tiny $\#$ of points} & 1002&1962 & 4002 &7842&16002\\ \hline\hline
\multicolumn{6}{|c|}{finite element method from all the point and mesh}\\ \hline
$\lambda=20$ & 0.0165 & 0.0085& 0.0042 & 0.0021 &0.0010\\ \hline
$\lambda=72$ & 0.0660 & 0.0342 & 0.0169 &0.0087 & 0.0043 \\
\hline\hline
\multicolumn{6}{|c|}{$\gamma = 100\%,\qquad  \ell=6$}\\ \hline
$\lambda=20$ & 0.0196 & 0.0084& 0.0039 & 0.0021 &0.0011\\ \hline
$\lambda=72$ & 0.1045 & 0.0599 & 0.0359 &0.0221 & 0.0139 \\
\hline\hline
\multicolumn{6}{|c|}{$\gamma = 90\%,\qquad  \ell= 10$}\\ \hline
$\lambda=20$ & 0.0494 & 0.0270 & 0.0137 & 0.0071 &0.0033\\ \hline
$\lambda=72$ & 0.2576 & 0.0792 & 0.0472 &0.0258 & 0.0143 \\
\hline\hline
\multicolumn{6}{|c|}{$\gamma = 85\%,\qquad  \ell= 20$}\\ \hline
$\lambda=20$ & 0.0647 & 0.0332 & 0.0154 & 0.0075 &  0.0040 \\ \hline
$\lambda=72$ & 0.4830 & 0.2020 & 0.0691 & 0.0369 & 0.0193 \\
\hline\hline
\multicolumn{6}{|c|}{$\gamma = 80\%,\qquad  \ell= 30$}\\ \hline
$\lambda=20$ & 0.0813 & 0.0397 & 0.0164& 0.0078 &0.0040\\ \hline
$\lambda=72$ & 0.6165 & 0.3387 & 0.0952&0.0483 & 0.0220 \\
\hline\hline
\multicolumn{6}{|c|}{3\% of global distance information ($\ell=6$ for MLS)}\\ \hline
$\lambda=20$ & 0.0155 & 0.0079& 0.0047 & 0.0030 &0.0021\\ \hline 
$\lambda=72$ & 0.0970 & 0.0550 & 0.0213 &0.0099 & 0.0069\\
\hline
\end{tabular}
}
\end{minipage}
\caption{Computation errors for LB eigenvalues $\lambda=20, 72$ on the unit sphere based on incomplete distance $D_{\gamma,\ell}$ without knowing exact KNN information. Left: convergence curves. Right: Relative errors for input distance matrices with different sample size.}
\label{fig:RelErrorsforLBeigenvaluesOnS2woKNN}
\end{figure}

\begin{figure}[h]
\centering
\begin{minipage}{0.38\linewidth}
\centering {\footnotesize $\lambda = 20$}\\
\includegraphics[width=1\linewidth]{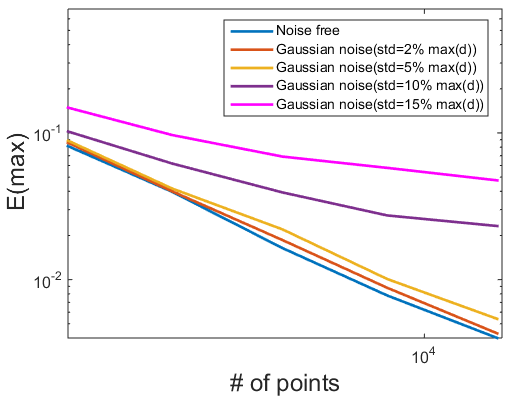}\\
\centering {\footnotesize $\lambda = 72$}\\
\includegraphics[width=1\linewidth]{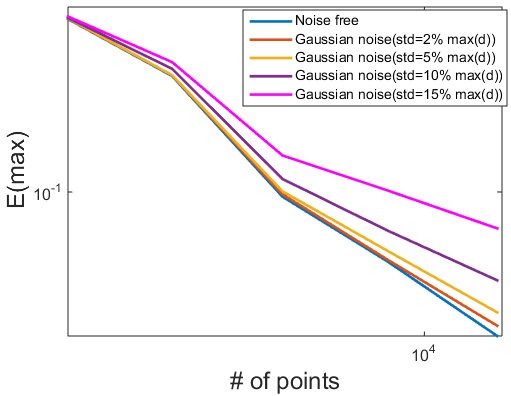}\\
\end{minipage}\hfill
\begin{minipage}{0.6\linewidth}
{\footnotesize
\begin{tabular}{|c|c|c|c|c|c|}
\hline
{\tiny $\#$ of points} & 1002&1962 & 4002 &7842&16002\\ \hline\hline
\multicolumn{6}{|c|}{Noise free,\quad  $\gamma = 100\%, \quad \ell=30$}\\ \hline
$\lambda=20$ & 0.0813 & 0.0397 & 0.0164& 0.0078 &0.0040\\ \hline
$\lambda=72$ & 0.6165 & 0.3387 & 0.0952&0.0483 & 0.0220 \\
\hline \hline
\multicolumn{6}{|c|}{ $\sigma= 2 \% \cdot d_{\max}, \quad \gamma = 80\%, \quad \ell=30$}\\ \hline
$\lambda=20$ & 0.0858 & 0.0400 & 0.0186 & 0.0088 &0.0043\\ \hline 
$\lambda=72$ & 0.6185 & 0.3410 & 0.0976 &0.0492 & 0.0245 \\ 
\hline\hline
\multicolumn{6}{|c|}{ $\sigma= 5 \% \cdot d_{\max}, \quad \gamma = 80\%, \quad \ell=30$}\\ \hline
$\lambda=20$ & 0.0886 & 0.0418 & 0.0220  & 0.0101 &0.0054\\ \hline 
$\lambda=72$ & 0.6218 & 0.3418 & 0.1007 &0.0540 & 0.0282 \\ 
\hline\hline
\multicolumn{6}{|c|}{ $\sigma= 10 \% \cdot d_{\max}, \quad \gamma = 80\%, \quad \ell=30$}\\ \hline
$\lambda=20$ & 0.1023 & 0.0619 & 0.0393 & 0.0274 &0.0232\\ \hline 
$\lambda=72$ & 0.6248 & 0.3653 & 0.1147 &0.0668 & 0.0395 \\ 
\hline\hline
\multicolumn{6}{|c|}{$\sigma= 15 \% \cdot d_{\max}, \quad \gamma = 80\%, \quad \ell=30$}\\ \hline
$\lambda=20$ & 0.1488 & 0.0968 & 0.0690 & 0.0578 &0.0475\\ \hline 
$\lambda=72$ & 0.6334 & 0.3920 & 0.1468 &0.1018 & 0.0682 \\ 
\hline
\end{tabular}
}
\end{minipage}
\caption{Computation errors for LB eigenvalues $\lambda=20, 72$ on the unit sphere based on incomplete corrupted distance $D_{\gamma,\ell} +  N(0,\sigma^2)$ without KNN information. Left: convergence curves. Right: Relative errors for input distance matrices with different sample size.}\label{fig:RelErrorsforLBeigenvaluesOnS2NoiseD}
\end{figure}

\begin{figure}[H]
\centering
\begin{tabular}{c@{\hspace{10pt}}c@{\hspace{10pt}}c@{\hspace{10pt}}c}
\includegraphics[width=0.2\linewidth]{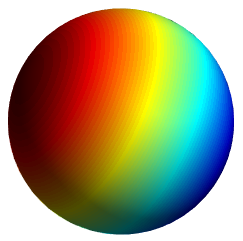}&
\includegraphics[width=0.2\linewidth]{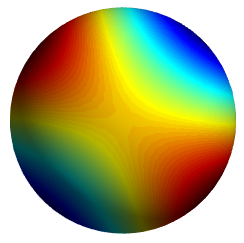}&
\includegraphics[width=0.2\linewidth]{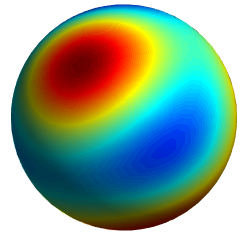}&
\includegraphics[width=0.2\linewidth]{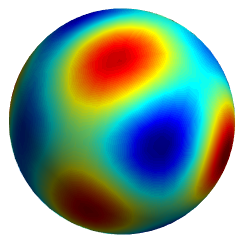}\\
\includegraphics[width=0.2\linewidth]{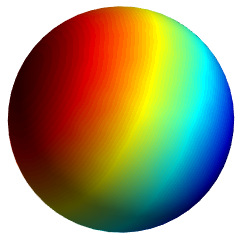}&
\includegraphics[width=0.2\linewidth]{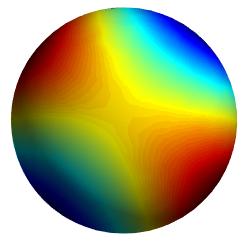}&
\includegraphics[width=0.2\linewidth]{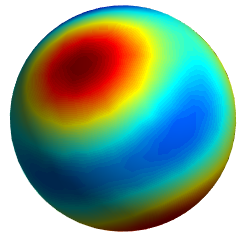}&
\includegraphics[width=0.2\linewidth]{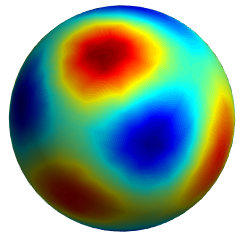}\\
\includegraphics[width=0.2\linewidth]{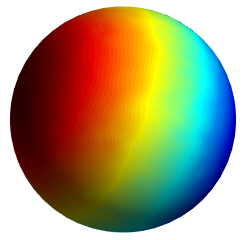}&
\includegraphics[width=0.2\linewidth]{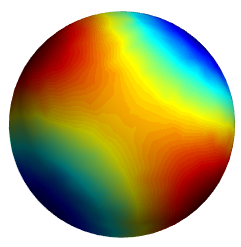}&
\includegraphics[width=0.2\linewidth]{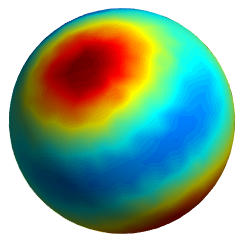}&
\includegraphics[width=0.2\linewidth]{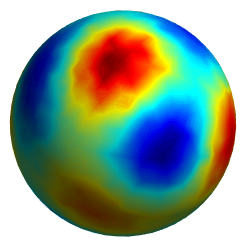}\\
\end{tabular}
\caption{LB eigenfunctions are color-coded on the unit sphere (1002 points) based on incomplete distance $D_{80\%,30}$ corrupted by Gaussian noise of different levels. Images from top to bottom represent the LB eigenfunctions from noise free distance, distance with Gaussian noise of $\sigma= 5 \% \cdot d_{\max} $, distance with Gaussian noise of $\sigma= 10 \% \cdot d_{\max} $, respectively. Images from left to right represent the LB eigenfunctions corresponding to $\lambda=2,6,12,20$, respectively. }
\label{fig:LBeigenS2NoiseD}
\end{figure}

\begin{figure}[H]
\centering
\begin{tabular}{c@{\hspace{10pt}}c@{\hspace{10pt}}c}
\includegraphics[width=0.24\linewidth]{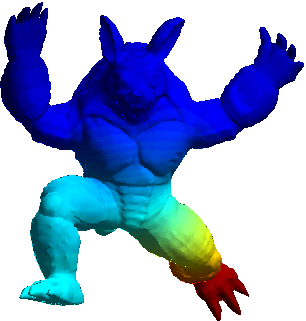}&
\includegraphics[width=0.24\linewidth]{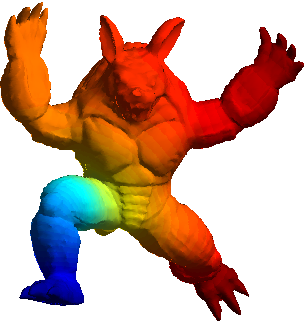}&
\includegraphics[width=0.24\linewidth]{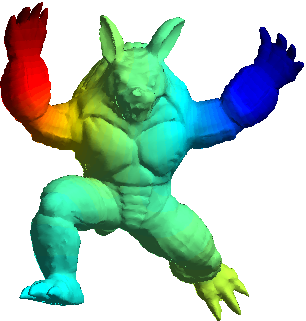}\\
\includegraphics[width=0.2\linewidth]{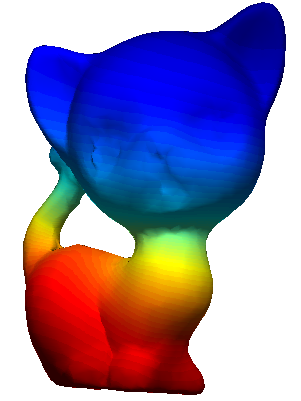}&
\includegraphics[width=0.2\linewidth]{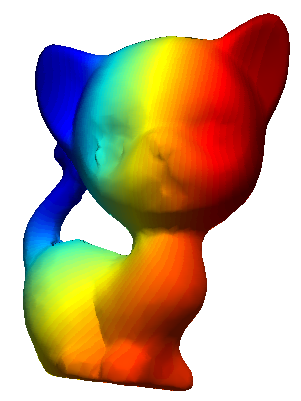}&
\includegraphics[width=0.2\linewidth]{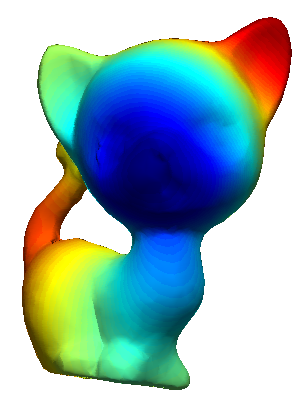}\\
\end{tabular}
\caption{LB eigenfunctions for the armadillo (16519 points) and the kitten (2884 points) surfaces based on incomplete distance $D_{50\%, 30}$. Top: the first three eigenfunctions are color-coded on the armadillo surface. Bottom: the first three eigenfunctions are color-coded on the kitten surfaces.}
\label{fig:LBeigenArmadilloKitten}
\end{figure}

In practice, the input condition might be even weaker than our previous experiment as the exact KNN information may not be directly obtained from incomplete local distances. Therefore, we also test our method for reconstructing local coordinates only relying on KNN information from incomplete distance information. In this scenario, we use the same setting as the previous test but require a larger rate of available distance information to have successfully local coordinate reconstruction. Our numerical results reported in Figure~\ref{fig:RelErrorsforLBeigenvaluesOnS2woKNN} also demonstrate nearly second order convergence of LB eigenvalue problem on the unit sphere. We further test the robustness of our method to distance with Gaussian perturbation. In this setting, we assume the input data is a incomplete distance $D_{80\%,30}$ corrupted by different levels of Gaussian noise with standard deviation $\sigma= 2\% d_{\max}, 5\% d_{\max}, 10\% d_{\max}, 15\% d_{\max}$ respectively, where $d_{\max} = \max_{i,j}\{D(i,j)\}$. The left table in Figure~\ref{fig:RelErrorsforLBeigenvaluesOnS2NoiseD} reports numerical accuracy for different levels of Gaussian noise. It is clear that our method still provides a reasonable good approximation of LB eigenvalues. Moreover, we also plot the corresponding LB eigenfunctions in Figure~\ref{fig:LBeigenS2NoiseD}, where LB eigenfunctions are color-coded on the unit sphere by setting red for positive values and blue for negative values. This figure illustrates consistent distribution patterns of LB eigenfunctions for noisy distance case. In addition, we also test our method of computing LB eigenvalue problems for incomplete distance information from point clouds sampled from more complicated manifolds such as armadillo (16519 points) and kitten (2884 points) surfaces. Figure~\ref{fig:LBeigenArmadilloKitten} reports serval LB eigenfunctions on armadillo and kitten surfaces based on incomplete distance matrix $D_{50\%,30}$.

\begin{figure}[h]
\centering
\begin{minipage}{0.4\linewidth}
\includegraphics[width=1\linewidth]{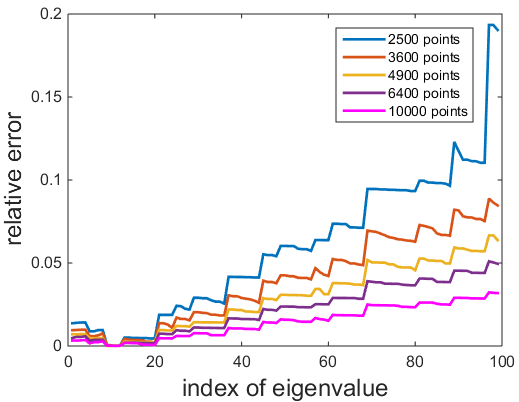}
\end{minipage}
\begin{minipage}{0.4\linewidth}
\includegraphics[width=1\linewidth]{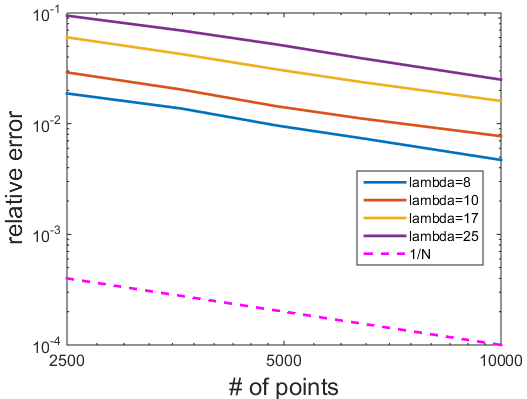}
\end{minipage}\hfill\\
\centering
\begin{minipage}{0.4\linewidth}
\includegraphics[width=1\linewidth]{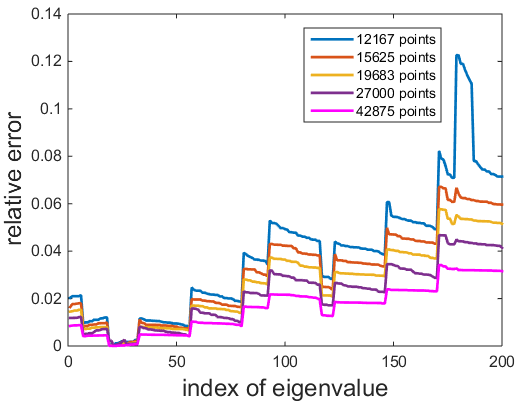}
\end{minipage}
\begin{minipage}{0.4\linewidth}
\includegraphics[width=1\linewidth]{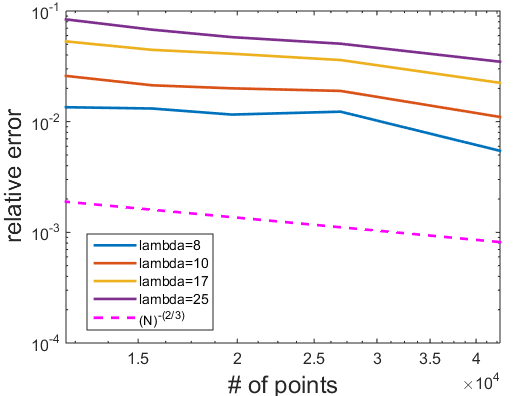}
\end{minipage}\hfill\\
\caption{Top: The first 100 LB eigenvalues for $T^2$ given by $D_{80\%,30}$ with different sample size $n = 2500, 3600, 4900, 6400, 10000$. Top left: Relative errors. Top right: Convergence curves. Bottom: The first 200 LB eigenvalues for $T^3$ given by $D_{85\%,30}$ with different sample size $n = 12167, 15625, 19683, 27000, 42875$. Bottom left: Relative errors. Bottom right: Convergence curves. }
\label{fig:LBeigsFlatTorus}
\end{figure}

As we discussed before, our methods can also be applied on manifolds with dimension more than 2 and co-dimension more than 1. Our third experiment tests the proposed method for computing LB eigenvalue problem on a 2 dimensional torus $T^2 = \{(\cos\theta_1,\sin\theta_1,\cos\theta_2,\sin\theta_2)\in \mathbb{R}^4~|~ \theta_1,\theta_2\in [0, 2\pi)\}$ in $\mathbb{R}^4$ and a 3 dimension torus $T^3 = \{(\cos\theta_1,\sin\theta_1,\cos\theta_2,\sin\theta_2, \cos\theta_3,\sin\theta_3)\in \mathbb{R}^6~|~ \theta_1, \theta_2, \theta_3\in [0, 2\pi)\}$ in $\mathbb{R}^6$. Solutions of the corresponding LB eigenvalue problem for both objects have closed forms, which enable us to conduct computation accuracy comparisons. For $T^2$, we choose $D_{80\%, 30}$ from the corresponding distance matrix obtained from different sample size on $T^2$. Similarly setting is considered for $T^3$ but using $D_{85\%,30}$ with slightly large portion of information as the rank for each local Gram matrix for $T^3$ is higher than the one for $T^2$. This requires more information of the distance matrix from the matrix completion theory. Nevertheless, we only use less than $1\%$ of the distance information for our experiments as $\ell = 30, n \geq 2500$. We remark that our computation does not assume the exact KNN information is available. This is the reason that we need higher portion of local information for accurate local coordinate reconstruction. In addition, we do not assume the dimension information is available in our computation. As we demonstrated in the first experiment, the value of $\gamma$ can be even lower if the exact KNN information is provided. We compute the first 100 LB eigenvalues for $T^2$ and $T^3$ using incomplete distance matrices with different sample size . Figure \ref{fig:LBeigsFlatTorus} reports the convergence and relative errors of our method. Both cases demonstrate approximately second order convergence.

\subsection{Solve Eikonal equation from incomplete distance}

\begin{figure}[h]
\centering
\begin{tabular}{cccc}
\includegraphics[width=0.249\linewidth]{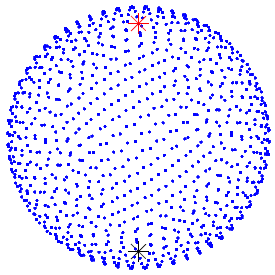}&
\includegraphics[width=0.249\linewidth]{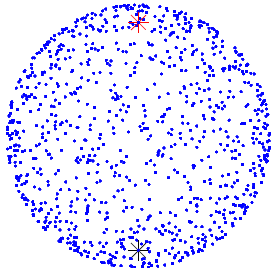}&
\includegraphics[width=0.19\linewidth]{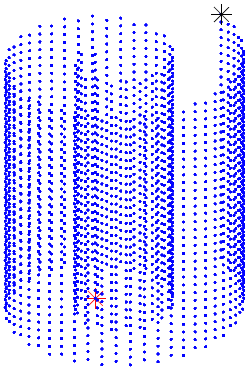}&
\includegraphics[width=0.19\linewidth]{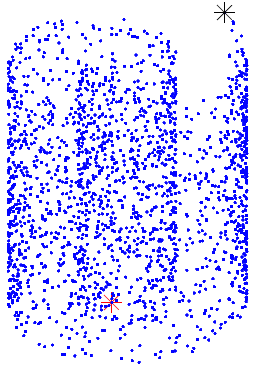}\\
\end{tabular}\\
\vspace{0.3cm}
\centering
{\footnotesize
\begin{tabular}{|cc||c|c|c|c|c|}
\hline
\multicolumn{2}{|c||}{\backslashbox{methods}{sample size}}& 1002&1962 & 4002 &7842&16002\\ \hline \hline
& &  \multicolumn{5}{|c|}{\textbf{Uniform sampling on $S^2$} }\\ \hline
\multicolumn{1}{ |c|  }{\multirow{2}{*}{Dijkstra} } & $E_{av}$ & 0.0348 & 0.0285 & 0.0248 & 0.0232 &0.0224\\
\cline{2-7}
\multicolumn{1}{ |c|  }{}     & $E_{se}$  & 0.008615 & 0.008606 & 0.008296 & 0.010642 &0.011501\\ \hline
\multicolumn{1}{ |c|  }{\multirow{2}{*}{Our method}}  & $E_{av}$ & 0.0113 & 0.0095 & 0.0080 & 0.0064 &0.0022\\
\cline{2-7}
\multicolumn{1}{ |c|  }{}  & $E_{se}$& 0.008100 & 0.005890 & 0.004110 & 0.002877 &0.002158\\ \hline\hline

&&\multicolumn{5}{|c|}{ \textbf{Non-Uniform sampling on $S^2$}}\\ \hline
\multicolumn{1}{ |c|  }{\multirow{2}{*}{Dijkstra} } & $E_{av}$& 0.0363 & 0.0344 & 0.0319 & 0.0305 &0.0294\\
\cline{2-7}
\multicolumn{1}{ |c|  }{}  & $E_{se}$& 0.011209 & 0.016090 & 0.018380 &0.016391  &0.019953\\ \hline 
\multicolumn{1}{ |c|  }{\multirow{2}{*}{Our method}}   & $E_{av}$& 0.0200 &0.0163  & 0.0141 & 0.0124 &0.0088\\
\cline{2-7}
\multicolumn{1}{ |c|  }{}  & $E_{se}$&0.012016& 0.008792& 0.003742  & 0.001736 &0.002765\\ \hline \hline

&&\multicolumn{5}{|c|}{ \textbf{Uniform sampling on a Swiss roll }}\\ \hline
\multicolumn{1}{ |c|  }{\multirow{2}{*}{Dijkstra} }  & $E_{av}$&0.0119  & 0.0156 & 0.0198&0.0200&0.0203\\
\cline{2-7}
 \multicolumn{1}{ |c|  }{}  & $E_{se}$& 0.013104 & 0.021242 & 0.024560 &0.024311&0.026004\\ \hline
\multicolumn{1}{ |c|  }{\multirow{2}{*}{Our method}} & $E_{av}$ & 0.0065 &  0.0044 & 0.0052&0.0033&0.0022\\
\cline{2-7}
\multicolumn{1}{ |c|  }{}   & $E_{se}$ & 0.003127 & 0.001637 & 0.001130&0.000783&0.000620\\ \hline\hline

&&\multicolumn{5}{|c|}{ \textbf{Non-Uniform sampling on a Swiss roll } }\\ \hline
\multicolumn{1}{ |c|  }{\multirow{2}{*}{Dijkstra} } & $E_{av}$ &0.0400  &  0.0105& 0.0114 &0.0151&0.0180\\
\cline{2-7}
\multicolumn{1}{ |c|  }{}   & $E_{se}$& 0.016612 &  0.015779&  0.014573&0.016587&0.018649\\ \hline
\multicolumn{1}{ |c|  }{\multirow{2}{*}{Our method}}  & $E_{av}$  & 0.0138 & 0.0053  & 0.0047&0.0029&0.0073\\
\cline{2-7}
\multicolumn{1}{ |c|  }{} & $E_{se}$  & 0.004754 & 0.005189 & 0.003087&0.005171&0.007246\\ \hline
\end{tabular}
}
\caption{Top left two images:  Incomplete distance data $D_{60\%,20}$ from 1002 points sampled uniformly and non-uniformly on the unit sphere with fixed north pole (red star) and south pole (blue star) for calculating geodesic distance. Top right two images: $D_{60\%,20}$ from 2048 points sampled uniformly and non-uniformly on the swiss role with a fixed starting point (red star) and an ending point (blue star) for calculating geodesic distance. Bottom table: Relative error of geodesic distances $E_{av}$ (averaging from starting to all the points) and $E_{se}$ (from the starting point to the ending point).}
\setlength{\abovecaptionskip}{10pt}
\label{fig:EikonalEqn}
\end{figure}

In this subsection, we test our method for solving a special hyperbolic equation, the Eikonal equation, on manifolds represented as incomplete distance information. In our experiments, assuming the KNN information is given as prior knowledge, we considered the incomplete local distance of uniform/non-uniform sampled point cloud of a unit sphere, in which the number of points vary from $1002$ to $16002$. In particular, we assume the point cloud always contains two points with coordinate $(0,0,1)$ and $(0,0,-1)$, namely, the north pole and the south pole  (See Figure \ref{fig:EikonalEqn}). To evaluate the quality of the reconstruction of distance map, we compare the relative error of the geodesic distance from the north pole to south pole (referred as $E_{se}$), and from the north pole to all other points (referred as $E_{av}$). The table in Figure \ref{fig:EikonalEqn} shows that our approach based on the fast marching method is much more accurate than the Dijkstra's method~\cite{dijkstra1959note}. Moreover, relative errors of our method decrease as the number of points increase. Furthermore, for the uniform sample cases, the relative error of geodesic distance from north pole to south pole has first order convergence with respect to the density of points using our method, while the Dijkstra's method does not have the same convergence property.

Similarly, we also compute the distance function on a Swiss roll: $$S = \{(t+0.1)\cos(t), (t+0.1)\sin(t), 8 \pi s)~|~t \in [0, 4 \pi], s \in [0,1] \}$$
 based on incomplete local distance from uniform/non-uniform sampling, in which two diagonal points are fixed and their position can be seen in Figure \ref{fig:EikonalEqn}. Using $D_{60\%,20}$ from different sample sizes vary from $1002$ to $16002$, the bottom Table in Figure \ref{fig:EikonalEqn} shows that our method can approximate the distance function with much less error than the Dijkstra's method, and approximate first convergence can also be observed from the case of the uniform sampled Swiss roll.

\subsection{Global coordinate reconstruction from using patch stitching}
\begin{figure}[h]
\centering
\begin{minipage}{0.34\linewidth}
\includegraphics[width=1\linewidth]{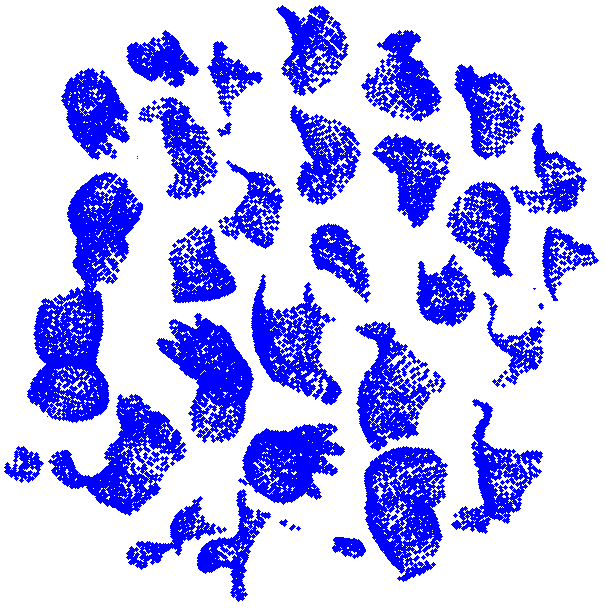}
\end{minipage}
\begin{minipage}{0.3\linewidth}
\includegraphics[width=1\linewidth]{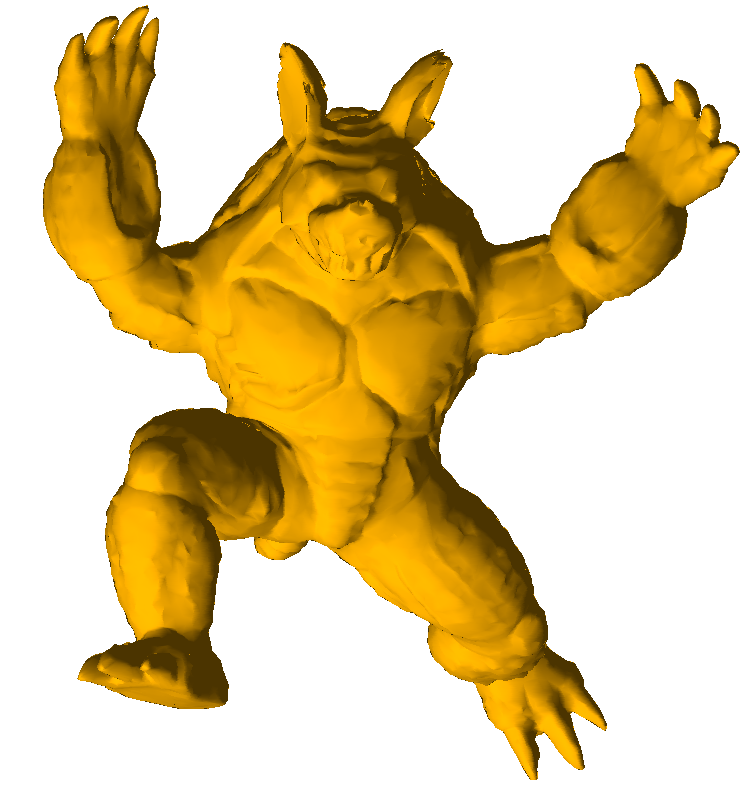}
\end{minipage}
\begin{minipage}{0.325\linewidth}
\includegraphics[width=1\linewidth]{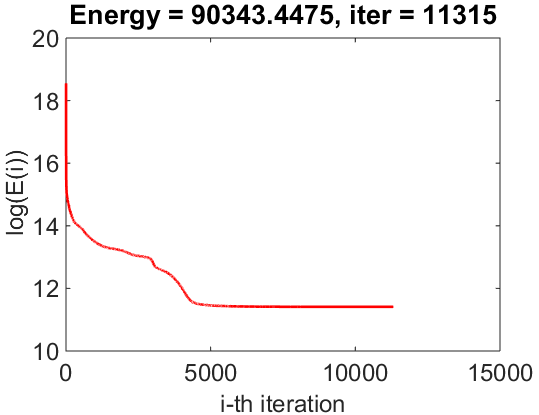}
\end{minipage}\\
\begin{minipage}{0.34\linewidth}
\includegraphics[width=1\linewidth]{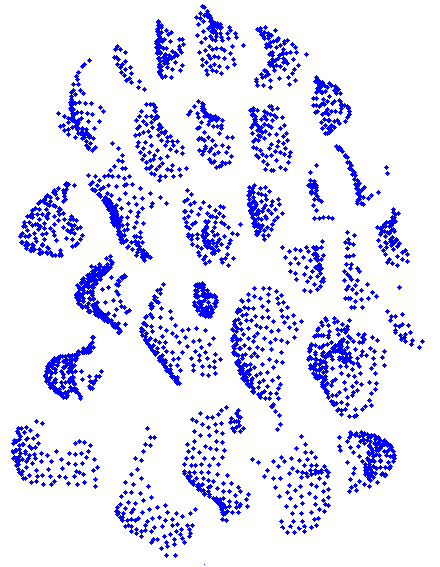}
\end{minipage}
\begin{minipage}{0.3\linewidth}
\includegraphics[width=.9\linewidth]{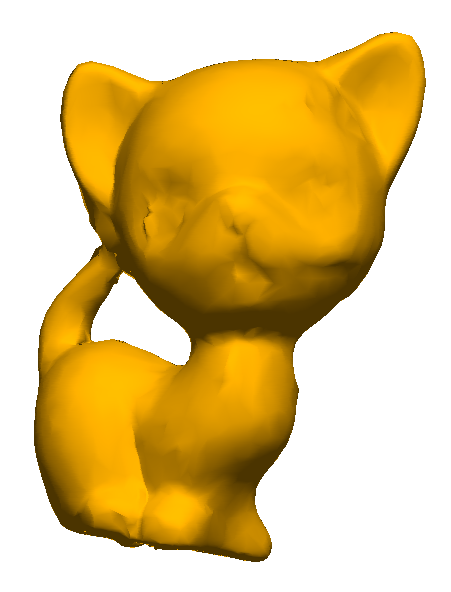}
\end{minipage}
\begin{minipage}{0.325\linewidth}
\includegraphics[width=1\linewidth]{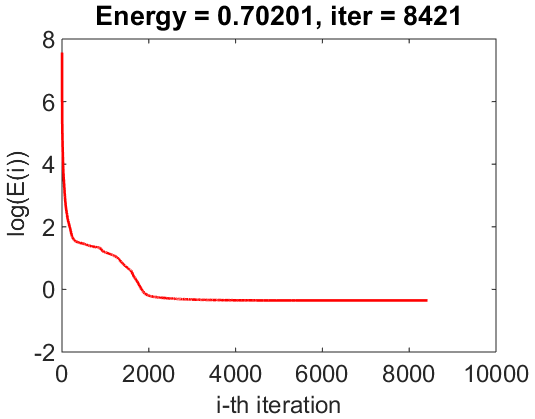}
\end{minipage}
\caption{Top: Local (left) and global (middle) coordinate reconstruction of the armadillo surface (16519 points, 28 patches) and the corresponding convergence curve of the global reconstruction objective value $\log \mathcal{E}$. Bottom: Local (left) and global (middle) coordinate reconstruction of the kitten surface( 2884 points, 30 patches) and the corresponding convergence curve of the global reconstruction objective value $\log \mathcal{E}$.}
\label{fig:stitching}
\end{figure}

In this experiment, we test our proposed global reconstruction model as an application of using global information from solving the LB eigenproblem based on incomplete distance. We assume the input manifold $\mathcal{M}\subset\RR^3$ sampled as a point cloud and has been separated as partially overlapped $L$ patches, namely, $\M = \bigcup_{j=1}^L \Omega_j$. We generate distance matrices $\{D_j\}_{j=1}^L$ from pairwise Euclidean distance on each patch and randomly choose $50\%$ information of each distance matrix. Based on these incomplete $\{D_j\}_{j=1}^L$, we first compute the first 100 LB eigenfunctions of $\M$ using the proposed method discussed in Section~\ref{sec:SolvePDEs}. After that, we reconstruct coordinates of each $\Omega_j$ from its incomplete distance using our matrix completion method. Finally, global reconstruction of $\M$ is obtained by the manifold stitching model discussed in Section~\ref{sec:MfdStitching}, which can handle the coordinate inconsistency from each local reconstruction. In Figure~\ref{fig:stitching}, we report numerical results for reconstructing an armadillo surface (with 16519 points, 28 patches) and a kitten surface (with 2884 points, 30 patches). The two left images in Figure~\ref{fig:stitching} illustrate coordinate reconstruction for each patch, the two middle images show the corresponding global reconstruction and the two right energy curves indicate the convergence of the proposed global reconstruction algorithm. It is clear to see that the proposed method can successfully reconstruct a manifold from a set of incomplete distance matrices from separated patches, while the global matrix completion model can not handle this case as the missing distance information is quite coherent. Moreover, compared to direct global coordinate reconstruction,  the method of local distance reconstruction and patch stitching saves more computational time. To validate this claim, we measured that the $50\%$ of local distance $\{D_j\}_{j=1}^L$ actually contains $3.59\%$ of full distance $D$ for armadillo surface and $3.45\%$ for kitten surface. Then we also implement the direct global coordinate reconstruction for these two surfaces using the same rate of random missing global distance $D$. Table \ref{fig:Time Comparison Local Global} shows that the computation time is much smaller using the stitching scheme especially for the armadillo surface with a large number of points.

\begin{table}
\centering
\begin{tabular}{|c|c|c|c|c|}
\hline
\multicolumn{2}{|c||}{\backslashbox{methods}{data}}& Armadillo& Kitten & Swiss roll \\
\hline
\multicolumn{2}{|c||}{Global reconstruction} & 35321.84 & 1315.50 & 622.65\\
\hline
\multicolumn{2}{|c||}{Stitching} & 760.18 & 138.76  & 199.75\\
\hline
\end{tabular}
\caption{Time comparisons (seconds) between global reconstruction and the stitching model.}\label{fig:Time Comparison Local Global}
\end{table}

\begin{figure}[h]
\centering
\begin{tabular}{c@{\hspace{10pt}}c}
\includegraphics[width=0.33\linewidth]{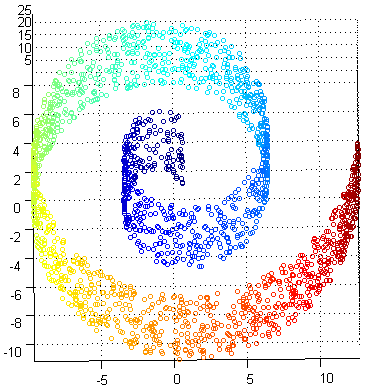}&
\includegraphics[width=0.4\linewidth]{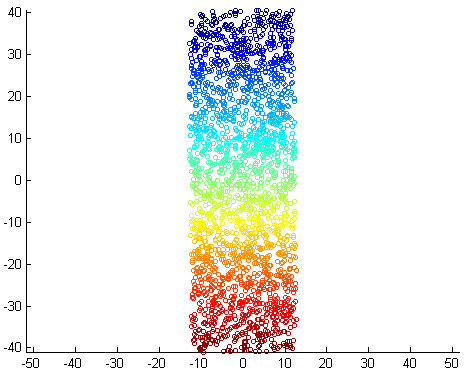}\\
\includegraphics[width=0.4\linewidth]{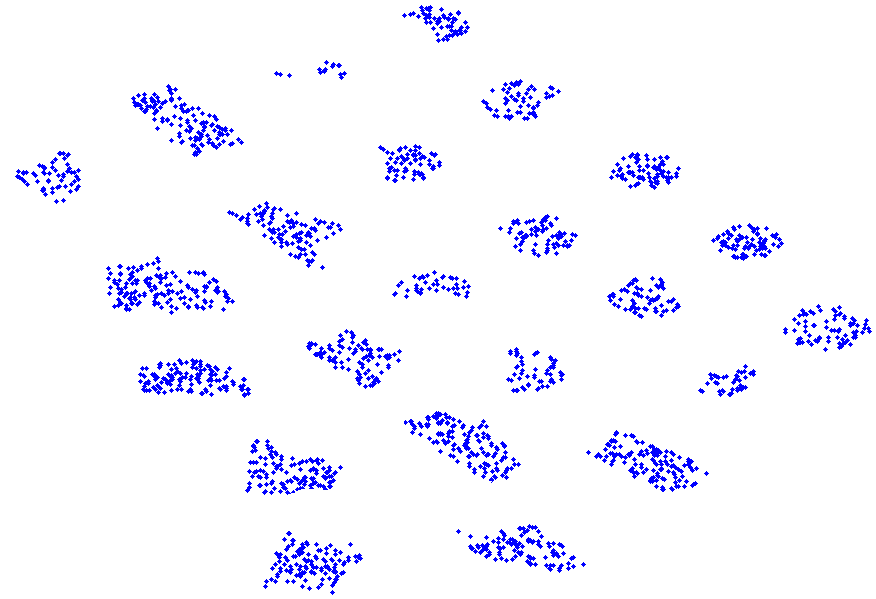}&
\includegraphics[width=0.4\linewidth]{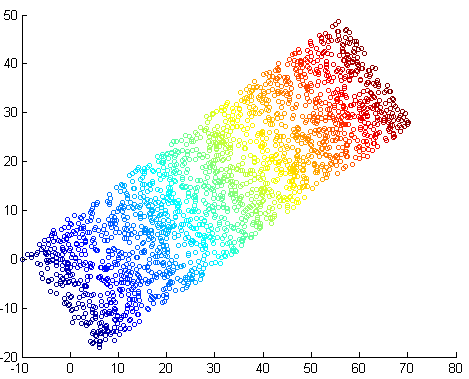}\\
\end{tabular}
\caption{Top: the Swiss roll surface (left) and its dimensional reduction result (right) from randomly $3.61\%$ of pair-wise geodesic distance. Bottom: local (left) and global (right) coordinates reconstruction of the Swiss roll from its  $80\%$ local geodesic distance.}
\label{fig:DimReduct_swissroll}
\end{figure}

We further conduct computation based on an incomplete geodesic distance matrix provided by a set of uniformly sampled points on the Swiss roll $S$ used in the previous example.  It is straightforward to check that the geodesic distance between any two points $((t_1+0.1)\cos(t_1), (t_1+0.1)\sin(t_1), 8 \pi s_1)$ and $((t_2+0.1)\cos(t_2), (t_2+0.1)\sin(t_2), 8 \pi s_2)$ is given by $ \sqrt{(8 \pi (s_1-s_2))^2 + (\mathcal{G}(t_2)-\mathcal{G}(t_1))^2} $, where $\mathcal{G}(t) = \frac{1}{2}[(t+0.1) \sqrt{1+(t+0.1)^2}+\log(|\sqrt{1+(t+0.1)^2}+(t+0.1)|)]$. In fact, this Swiss roll $S$ is isometrical to a 2D flat domain $[\mathcal{G}(0), \mathcal{G} (4 \pi)] \times [0, 8 \pi]$. We uniformly sample 2048 points on $S$ and construct a distance matrix $D$ using pairwise geodesic distance. We test two ways of reconstruct the Swiss roll based on incomplete information of $D$. Firstly, similar as patch stitching simulations, we assume the Swiss roll is separated as $23$ patches with partial overlap. Then we chose only $50\%$ of local geodesic distance $\{D_j\}_{j=1}^{23}$, which is essentially $3.61\%$ coherent sampling of totally geodesic distance $D$. The bottom row of figure~\ref{fig:DimReduct_swissroll} shows that the global stitched coordinate is a 2D rectangle, which is identical to the theoretical ground truth of the dimension reduced Swiss roll. Secondly, we assume $3.61\%$ random sampled global geodesic distance is available and apply our coordinate reconstruction algorithm to have the global coordinate reconstruction for the Swiss roll. Top right image in figure~\ref{fig:DimReduct_swissroll} shows that direct global coordinate reconstruction from incoherent random missing distance is also applicable. Table \ref{fig:Time Comparison Local Global} shows that with same sampling rate from the whole distance matrix, the direct global coordinate reconstruction takes much more computational time. The reason can be explained similarly as computation complexity discussion in subsections~\ref{subsec:ReviewMLSLM} and ~\ref{subsec:SovelLBeigs}. The step of solving PDEs from local distance has time complexity $\mathcal{O}(\min(l^2mn,l^3n))$ while direct global coordinate reconstruction has time complexity $\mathcal{O}(mn^2)$ which is larger if total number of points $n$ is large enough. Although the eigen-decomposition of Laplace-Beltrami matrix also involves global computation but it needs to be performed only once while eigenvalue thresholding is required in every iteration for direct global coordinate reconstruction. Therefore, it is not uprising that the stitching scheme saves a lot of time. 
We remark that our method of using geodesic distance to reconstructing manifolds is highly related to an important topic called nonlinear dimension reduction in machine learning and statistics community~\cite{roweis2000nonlinear}. In our future work, we will investigate more along this direction based on our methods.

\section{Conclusion}
\label{sec:conclusion}
In this paper, we proposed a framework for discretizing PDEs on manifolds represented as incomplete inter-point distance information. Our approach is to conduct PDE discretization point-wisely from a local coordinate reconstruction. This strategy successfully avoids a time-consuming step of global coordinate reconstruction and leads to a method with complexity linearly scaling to the number of sample size. Our local reconstruction model is inspired from the recent advances of low-rank matrix completion theory, which guarantee accuracy local coordinate reconstruction by only requiring a very small portion of distance information. Our method can be viewed as natural extensions of the moving least square method and the local method method for solving PDEs on point clouds to manifold-structured data represented as incomplete inter-point distance information. As an application of using solutions of PDEs, we propose a new manifold reconstruction model by stitching local patches on the spectrum domain. Intense numerical experiments indicate effectiveness and efficiency of our methods for solving LB eigenvalue problem, Eikonal equations and manifold reconstruction using stitching methods based on global information from the LB eigensystem.

\section{Acknowledgement}
We would like to thank Prof. Xiangxiaong Zhang's suggestions about optimizing distance matrix directly as we discussed in Remark~\ref{rk:LowrankD}.

\bibliographystyle{siam}             
\bibliography{SolvePDEfromDistance}

\begin{thebibliography}{10}

\bibitem{alfakih1999solving}
{\sc A.~Y. Alfakih, A.~Khandani, and H.~Wolkowicz}, {\em Solving {E}uclidean
  distance matrix completion problems via semidefinite programming},
  Computational optimization and applications, 12 (1999), pp.~13--30.

\bibitem{axelsson1999processing}
{\sc P.~Axelsson}, {\em Processing of laser scanner data---algorithms and
  applications}, ISPRS Journal of Photogrammetry and Remote Sensing, 54 (1999),
  pp.~138--147.

\bibitem{Belkin:ML2004}
{\sc M.~Belkin and P.~Niyogi}, {\em Semi-supervised learning on {R}iemannian
  manifolds}, Machine Learning, 56 (2004), pp.~209--239.

\bibitem{Belkin:09clp}
{\sc M.~Belkin, J.~Sun, and Y.~Wang}, {\em Constructing {L}aplace operator from
  point clouds in $\mathbb{R}^d$}, in Proceedings of the Twentieth Annual
  ACM-SIAM Symposium on Discrete Algorithms, Philadelphia, PA, USA, 2009,
  pp.~1031--1040.

\bibitem{Berard:1994}
{\sc P.~B\'erard, G.~Besson, and S.~Gallot}, {\em Embedding {R}iemannian
  manifolds by their heat kernel}, Geometric \& Functional Analysis, 4 (1994),
  pp.~373--398.

\bibitem{berger1999reconstructing}
{\sc B.~Berger, J.~Kleinberg, and T.~Leighton}, {\em Reconstructing a
  three-dimensional model with arbitrary errors}, Journal of the ACM (JACM), 46
  (1999), pp.~212--235.

\bibitem{Bertalmio:2002}
{\sc M.~Bertalmio, L.-T. Cheng, S.~Osher, and G.~Sapiro}, {\em Variational
  problems and partial differential equations on implicit surfaces}, Journal of
  Computational Physics, 174 (2002), pp.~759--780.

\bibitem{Bertalmio:2000}
{\sc M.~Bertalmio, G.~Sapiro, L.-T. Cheng, and S.~Osher}, {\em A framework for
  solving surface partial differential equations for computer graphics
  applications}, UCLA CAM Report (00-43),  (2000).

\bibitem{biswas2006semidefinite}
{\sc P.~Biswas, T.-C. Liang, K.-C. Toh, Y.~Ye, and T.-C. Wang}, {\em
  Semidefinite programming approaches for sensor network localization with
  noisy distance measurements}, IEEE transactions on automation science and
  engineering, 3 (2006), pp.~360--371.

\bibitem{borg2005modern}
{\sc I.~Borg and P.~J. Groenen}, {\em Modern multidimensional scaling: Theory
  and applications}, Springer Science \& Business Media, 2005.

\bibitem{Brandman:2008JSC}
{\sc J.~Brandman}, {\em A level-set method for computing the eigenvalues of
  elliptic operators defined on compact hypersurfaces}, Journal of Scientific
  Computing, 37 (2008), pp.~282--315.

\bibitem{Bronstein:2010CVPR}
{\sc M.~M. Bronstein and I.~Kokkinos}, {\em Scale-invariant heat kernel
  signatures for non-rigid shape recognition}, IEEE Conference on Computer
  Vision and Pattern Recognition (CVPR),  (2010), pp.~1704--1711.

\bibitem{CandesRe2008}
{\sc E.~Cand{\`e}s and B.~Recht}, {\em Exact matrix completion via convex
  optimization.}, Found. of Comput. Math., 9 (2008), pp.~717--772.

\bibitem{chaudhury2015global}
{\sc K.~N. Chaudhury, Y.~Khoo, and A.~Singer}, {\em Global registration of
  multiple point clouds using semidefinite programming}, SIAM Journal on
  Optimization, 25 (2015), pp.~468--501.

\bibitem{Chavel:1984}
{\sc I.~Chavel}, {\em Eigenvalues in Riemannian geometry}, Academic press. INC,
  1984.

\bibitem{coifman2006diffusion}
{\sc R.~R. Coifman and S.~Lafon}, {\em Diffusion maps}, Applied and
  computational harmonic analysis, 21 (2006), pp.~5--30.

\bibitem{crippen1988distance}
{\sc G.~M. Crippen, T.~F. Havel, et~al.}, {\em Distance geometry and molecular
  conformation}, vol.~74, Research Studies Press Taunton, UK, 1988.

\bibitem{cucuringu2012sensor}
{\sc M.~Cucuringu, Y.~Lipman, and A.~Singer}, {\em Sensor network localization
  by eigenvector synchronization over the euclidean group}, ACM Transactions on
  Sensor Networks (TOSN), 8 (2012), p.~19.

\bibitem{dijkstra1959note}
{\sc E.~W. Dijkstra}, {\em A note on two problems in connexion with graphs},
  Numerische mathematik, 1 (1959), pp.~269--271.

\bibitem{dziuk2013finite}
{\sc G.~Dziuk and C.~M. Elliott}, {\em Finite element methods for surface
  {PDE}s}, Acta Numerica, 22 (2013), pp.~289--396.

\bibitem{faugeras1986representation}
{\sc O.~Faugeras and M.~Hebert}, {\em The representation, recognition, and
  locating of 3-d objects}, The international journal of robotics research, 5
  (1986), pp.~27--52.

\bibitem{ji2004sensor}
{\sc X.~Ji and H.~Zha}, {\em Sensor positioning in wireless ad-hoc sensor
  networks using multidimensional scaling}, in INFOCOM 2004. Twenty-third
  AnnualJoint Conference of the IEEE Computer and Communications Societies,
  vol.~4, IEEE, 2004, pp.~2652--2661.

\bibitem{jolliffe2002principal}
{\sc I.~Jolliffe}, {\em Principal component analysis}, Wiley Online Library,
  2002.

\bibitem{Peter:08}
{\sc P.~W. Jones, M.~Maggioni, and R.~Schul}, {\em Manifold parametrizations by
  eigenfunctions of the {L}aplacian and heat kernels}, Proceedings of the
  National Academy of Sciences, 105 (2008), pp.~1803--1808.

\bibitem{kimmel1998computing}
{\sc R.~Kimmel and J.~A. Sethian}, {\em Computing geodesic paths on manifolds},
  Proceedings of the National Academy of Sciences, 95 (1998), pp.~8431--8435.

\bibitem{kruskal1978multidimensional}
{\sc J.~B. Kruskal and M.~Wish}, {\em Multidimensional scaling}, vol.~11, Sage,
  1978.

\bibitem{lai2011framework}
{\sc R.~Lai and T.~F. Chan}, {\em A framework for intrinsic image processing on
  surfaces}, Computer vision and image understanding, 115 (2011),
  pp.~1647--1661.

\bibitem{lai2017Exact}
{\sc R.~Lai, J.~Li, and A.~Tasissa}, {\em Exact reconstruction of distance
  geometry problem using low-rank matrix completion}, Preprint,  (2017).

\bibitem{lai2013local}
{\sc R.~Lai, J.~Liang, and H.~Zhao}, {\em A local mesh method for solving
  {PDE}s on point clouds.}, Inverse Problems \& Imaging, 7 (2013).

\bibitem{Lai:2010CVPR}
{\sc R.~Lai, Y.~Shi, K.~Scheibel, S.~Fears, R.~Woods, A.~W. Toga, and T.~F.
  Chan}, {\em Metric-induced optimal embedding for intrinsic 3{D} shape
  analysis}, Computer Vision and Pattern Recognition (CVPR),  (2010),
  pp.~2871--2878.

\bibitem{lai2017multi}
{\sc R.~Lai and H.~Zhao}, {\em Multiscale nonrigid point cloud registration
  using robust sliced-wasserstein distance via laplace-beltrami eigenmap}, SIAM
  Journal on Imaging Sciences, 10 (2017), pp.~449--483.

\bibitem{Levy:2006IEEECSMA}
{\sc B.~Levy}, {\em Laplace-{B}eltrami eigenfunctions: Towards an algorithm
  that understands geometry}, IEEE International Conference on Shape Modeling
  and Applications, invited talk,  (2006).

\bibitem{Liang:CVPR2012}
{\sc J.~Liang, R.~Lai, T.~Wong, and H.~Zhao}, {\em Geometric understanding of
  point clouds using {L}aplace-{B}eltrami operator}, CVPR,  (2012).

\bibitem{liang2013solving}
{\sc J.~Liang and H.~Zhao}, {\em Solving partial differential equations on
  point clouds}, SIAM Journal on Scientific Computing, 35 (2013),
  pp.~A1461--A1486.

\bibitem{lui2008variational}
{\sc L.~M. Lui, X.~Gu, T.~F. Chan, S.~Yau, et~al.}, {\em Variational method on
  riemann surfaces using conformal parameterization and its applications to
  image processing}, Methods and Applications of Analysis, 15 (2008),
  pp.~513--538.

\bibitem{Meyer:2003}
{\sc M.~Meyer, M.~Desbrun, P.~Schr\"oder, and A.~H. Barr}, {\em Discrete
  differential-geometry operators for triangulated 2-manifolds}, Visualization
  and Mathematics III. (H.C. Hege and K. Polthier, Ed.) Springer Verlag,
  (2003), pp.~35--57.

\bibitem{mucherino2012distance}
{\sc A.~Mucherino, C.~Lavor, L.~Liberti, and N.~Maculan}, {\em Distance
  geometry: theory, methods, and applications}, Springer Science \& Business
  Media, 2012.

\bibitem{Osher:88}
{\sc S.~Osher and J.~Sethian.}, {\em Fronts propagation with
  curvature-dependent speed: algorithms based on {H}amilton-{J}acobi
  formulations.}, Journal of computational physics, 79 (1988), pp.~12--49.

\bibitem{pinkall1993computing}
{\sc U.~Pinkall and K.~Polthier}, {\em Computing discrete minimal surfaces and
  their conjugates}, Experimental mathematics, 2 (1993), pp.~15--36.

\bibitem{recht2010guaranteed}
{\sc B.~Recht, M.~Fazel, and P.~A. Parrilo}, {\em Guaranteed minimum-rank
  solutions of linear matrix equations via nuclear norm minimization}, SIAM
  review, 52 (2010), pp.~471--501.

\bibitem{Reuter:06}
{\sc M.~Reuter, F.~Wolter, and N.~Peinecke}, {\em Laplace-{B}eltrami spectra as
  {S}hape-{DNA} of surfaces and solids}, Computer-Aided Design, 38 (2006),
  pp.~342--366.

\bibitem{roweis2000nonlinear}
{\sc S.~T. Roweis and L.~K. Saul}, {\em Nonlinear dimensionality reduction by
  locally linear embedding}, Science, 290 (2000), pp.~2323--2326.

\bibitem{saul2003think}
{\sc L.~K. Saul and S.~T. Roweis}, {\em Think globally, fit locally:
  unsupervised learning of low dimensional manifolds}, Journal of Machine
  Learning Research, 4 (2003), pp.~119--155.

\bibitem{scott2011sage}
{\sc J.~Scott and P.~J. Carrington}, {\em The SAGE handbook of social network
  analysis}, SAGE publications, 2011.

\bibitem{sethian1996fast}
{\sc J.~A. Sethian}, {\em {A fast marching level set method for monotonically
  advancing fronts}}, Proceedings of the National Academy of Sciences, 93
  (1996), pp.~1591--1595.

\bibitem{Shi:08a}
{\sc Y.~Shi, R.~Lai, S.~Krishna, N.~Sicotte, I.~Dinov, and A.~W. Toga}, {\em
  Anisotropic {L}aplace-{B}eltrami eigenmaps: bridging {R}eeb graphs and
  skeletons}, in Computer Vision and Pattern Recognition Workshops, 2008,
  pp.~1--7.

\bibitem{singer2008remark}
{\sc A.~Singer}, {\em A remark on global positioning from local distances},
  Proceedings of the National Academy of Sciences, 105 (2008), pp.~9507--9511.

\bibitem{spira2007geometric}
{\sc A.~Spira and R.~Kimmel}, {\em Geometric curve flows on parametric
  manifolds}, Journal of Computational Physics, 223 (2007), pp.~235--249.

\bibitem{stam2003flows}
{\sc J.~Stam}, {\em Flows on surfaces of arbitrary topology}, ACM Transactions
  On Graphics (TOG), 22 (2003), pp.~724--731.

\bibitem{starck1998image}
{\sc J.-L. Starck, F.~Murtagh, and A.~Bijaoui}, {\em Image processing and data
  analysis: the multiscale approach}, Cambridge University Press, 1998.

\bibitem{Sun:2009SGP}
{\sc J.~Sun, M.~Ovsjanikov, and L.~Guibas}, {\em A concise and provabley
  informative multi-scale signature baded on heat diffusion}, Eurographics
  Symposium on Geometry Processing,  (2009).

\bibitem{taubin2000geometric}
{\sc G.~Taubin et~al.}, {\em Geometric signal processing on polygonal meshes},
  Eurographics State of the Art Reports, 4 (2000), pp.~81--96.

\bibitem{SDO}
{\sc M.~J. Todd}, {\em Semidefinite optimization}, Acta Numerica, 10 (2001),
  pp.~515--560.

\bibitem{tseng2001convergence}
{\sc P.~Tseng}, {\em Convergence of a block coordinate descent method for
  nondifferentiable minimization}, Journal of optimization theory and
  applications, 109 (2001), pp.~475--494.

\bibitem{Vallet:2008CGF}
{\sc B.~Vallet and B.~Levy}, {\em Spectral geometry processing with manifold
  harmonics}, Computer Graphics Forum (Proceedings Eurographics),  (2008).

\bibitem{SDP}
{\sc L.~Vandenberghe and S.~Boyd}, {\em Semidefinite programming}, SIAM Review,
  38 (1996), pp.~49--95.

\bibitem{wang2007brain}
{\sc Y.~Wang, L.~M. Lui, .~Gu, K.~M. Hayashi, T.~F. Chan, A.~W. Toga, P.~M.
  Thompson, and S.~Yau}, {\em Brain surface conformal parameterization using
  riemann surface structure}, IEEE transactions on medical imaging, 26 (2007),
  pp.~853--865.

\bibitem{waterman1995introduction}
{\sc M.~Waterman}, {\em Introduction to computational biology: maps, sequences
  and genomes}, CRC Press, 1995.

\bibitem{weinberger2004learning}
{\sc K.~Q. Weinberger, F.~Sha, and L.~K. Saul}, {\em Learning a kernel matrix
  for nonlinear dimensionality reduction}, in Proceedings of the twenty-first
  international conference on Machine learning, ACM, 2004, p.~106.

\bibitem{AltSDP-WenGoldfarbYin-2009}
{\sc Z.~Wen, D.~Goldfarb, and W.~Yin}, {\em Alternating direction augmented
  lagrangian methods for semidefinite programming}, Mathematical Programming
  Computation, 2 (2010), pp.~203--230.

\bibitem{wen2013feasible}
{\sc Z.~Wen and W.~Yin}, {\em A feasible method for optimization with
  orthogonality constraints}, Mathematical Programming, 142 (2013),
  pp.~397--434.

\bibitem{xu2004convergent}
{\sc G.~Xu}, {\em Convergent discrete {L}aplace-{B}eltrami operators over
  triangular surfaces}, in Geometric Modeling and Processing, 2004.
  Proceedings, IEEE, 2004, pp.~195--204.

\bibitem{zhao2005fast}
{\sc H.~Zhao}, {\em A fast sweeping method for eikonal equations}, Mathematics
  of computation, 74 (2005), pp.~603--628.

\end{thebibliography}
\end{document}